\documentclass[10pt,a4paper]{article}
\usepackage{graphicx}
\usepackage{dcolumn}
\usepackage{bm}
\usepackage{graphics}
\usepackage[centerlast]{caption2}
\usepackage{epsfig}
\usepackage{color}
\usepackage{amssymb}
\usepackage{amsthm}
\usepackage{lineno}
\usepackage{mathrsfs}
\usepackage{amsmath}
\usepackage{amsfonts}
\usepackage{bm,cite}
\numberwithin{equation}{section}
\usepackage{tabularx}
\usepackage[top=2cm, bottom=2cm, left=3cm, right=3cm] {geometry}
\newtheorem{thm}{Theorem}[section]
\newtheorem{cor}[thm]{Corollary}
\newtheorem{lem}[thm]{Lemma}
\newtheorem{prop}[thm]{Proposition}
\theoremstyle{definition}

\theoremstyle{definition}

\theoremstyle{definition}

\def\R{{\mathbb{R}}}
\usepackage{lineno}

\begin{document}
	
\title{On Solution Uniqueness and Robust Recovery for Sparse Regularization with a Gauge:	from Dual Point of View}
\author{Jiahuan He\\School of Science,
	Harbin University of Science and Technology,\\
	Harbin, 150080, P.R. China\\
	Chao Kan and Wen Song\footnote{Corresponding author, E-mail address: wsong@hrbnu.edu.cn}\\School of Mathematical and Sciences, Harbin Normal University,\\ Harbin, 150025, P.R. China }
\date{}
\maketitle
\par\noindent
	 \textrm{\textbf{Abstract.}}
		In this paper, we focus  on the exploration of solution uniqueness, sharpness, and robust recovery in sparse regularization with a gauge $J$. 
	Based on the criteria for the uniqueness of Lagrange multipliers in the dual problem, we give a characterization of the unique solution via the so-called radial cone. We establish characterizations of the isolated calmness (i.e., uniqueness of solutions combined with the calmness properties) of two types of the solution mappings and prove that they are equivalent under the assumption of metric subregularity of the subdifferentials for regularizers. Furthermore, we present sufficient and necessary conditions for a sharp solution, and show that these conditions guarantee robust recovery with a linear rate and imply local upper Lipschitzian properties of the solution mapping. Some applications of the polyhedral regularizer case such as sparse analysis regularization, weighted sorted $\ell_1$-norm sparse regularization, and non-polyhedral regularizer cases such as nuclear norm optimization, conic gauge optimization, and semidefinite conic gauge optimization are given.\vspace*{2mm}
	\vspace*{2mm}
	
	\renewcommand{\thefootnote}{\fnsymbol{footnote}}
	\footnote[0]{	The research of the second author was supported by the Foundation of Heilongjiang Department
		of Education (UNPYSCT-2018181). The research of the third author was supported by the
		National Natural Sciences Grant of China (No. 11871182).}
	\section{Introduction}
	Consider a matrix $A \in \R^{m\times n}$ with $m \ll n$ and $b_0 \in \R^m$, and define the underdetermined linear system of
	equations
	\begin{equation}  \label{lin-equation} 
		Ax= b_0.
	\end{equation}
	This system has infinitely many solutions, which we denote as the solution set
	$\mathcal{X}:=\{ x \in \R^n \colon Ax= b_0\}$.
	To obtain the desired solution $x_0$, it is typical to address the following optimization problem:
	$$ \min_{x} J(x) \;\; \mbox{ s.t. } \;\; Ax=b_0, \leqno(P_J)$$
	where $J$ is known as a regularizer that is bounded from below. A straightforward and intuitive way to measure the sparsity of a solution is to solve the optimization problem with $J(x) = \|x\|_{0}$, that is
	$$ \min_{x} \|x\|_0  \;\; \mbox{ s.t. } \;\; Ax=b_0, \leqno(P_0) $$
	which is a combinatorial optimization problem with a prohibitive complexity if solved by enumeration.
	Alternatively, consider solving the $\ell_1$-norm regularization problem:
	$$ \min_{x} \|x\|_1  \;\; \mbox{ s.t. } \;\;  Ax = b_0. \leqno(P_1)$$
	This problem, known as basis pursuit, is an example of convexification of $(P_0)$. Other regularizers used in the literature include the weighted sorted $\ell_1$-norm, the weighted $\ell_1$ norm,the analysis sparsity seminorm, the $\ell_1/\ell_2$ norm, and the nuclear norm.

	When the observation is subject to errors, the system (\ref{lin-equation}) is modified to
	\begin{equation} \label{pertu-equation}
		b = Ax_0 + \omega,
	\end{equation}
	where $\omega$ accounts either for noise and/or modeling errors. Throughout, $\omega$ will be assumed deterministic with $\|\omega\|\le
	\delta$, with $\delta > 0$ known.
	
	To recover $x_0$ from (\ref{pertu-equation}), one can solve the noise-aware form (Mozorov regularization)
	\begin{equation} \label{robu-prob}
		\min_{x} J(x) \;\; \mbox{ s.t. } \;\;  \| Ax-b \| \le  \delta,
	\end{equation}
	or Tikhonov regularization
	\begin{equation} \label{regu-prob}
		\min_{x}  \frac{1}{2} \|Ax - b \|^2+ \mu J(x)
	\end{equation}
	with regularization parameter $\mu > 0$. It is well-known that problems (\ref{robu-prob}) and (\ref{regu-prob}) are formally equivalent in the sense that there exists a bijection between $\delta$ and $\mu$  such that both problems share the same set of solutions under certain additional conditions (see \cite{GHS,IVT}).
	
	In the fields of compressive sensing, signal and image processing, statistics, machine learning, inverse problems, and related areas, we are concerned with whether the solution is unique and the recovery is stable. We say that the recovery is robust when solutions to (\ref{robu-prob}) and (\ref{regu-prob}) are closely approximate the original signal $x_0$ as long as $\mu$ is appropriately chosen (as a function of the noise level $\delta$). Furthermore, linear rate robust recovery indicates that the error is a linear function of $\delta$.
	
	Several sufficient conditions have been suggested to ensure solution uniqueness and successful recovery. There are some uniform sufficient conditions guaranteeing solution uniqueness to problem $(P_0)$ or $(P_1)$ and recovery of all $s$-sparse solutions (i.e., solution $x$ to (\ref{lin-equation}) with $\|x\|_0 \le s$) as well as robust recovery  with a linear rate, such as the restricted isometry principle, the mutual coherence condition, and the null-space property (see \cite{Candes,CRT,FR}). The columns of $A$ are in general position ensuring a strong type of solution uniqueness to problem $(P_1)$ (i.e., uniqueness valid for all $b \in \R^m$) (see \cite{Dossal}) and problem (\ref{regu-prob}) for all $b\in \R^m$ and $\mu >0$ (see \cite{Tibshirani}). Recently, this condition was relaxed in \cite{ES} 	to a geometric criterion that is both sufficient and necessary for strong uniqueness in Lasso, and it was extended to a polyhedral norm in \cite{STW}.
	
	Fuchs \cite{Fuchs} provided a non-uniform sufficient condition for solution uniqueness to $(P_1)$, which is a dual type of condition. Later, its necessity was established in \cite{Dossal,GHS,ZYC}. Solution uniqueness has also been studied for more general regularizer functions, such as the polyhedral gauge \cite{Gilbert}, convex piecewise affine function \cite{MS}, and $\ell_1$-analysis regularizer \cite{FNT,ZYY}. Recently, \cite{LPB} gave the characterization of the unique solution to problem $(P_J)$ via the so-called radial cone and also two kinds of sufficient conditions for solution uniqueness when $J$ is the nuclear norm. Based on the subdifferential decomposability, \cite[Corollary 1]{VGFP} presented a sufficient condition for the solution uniqueness to problem $(P_J)$ whenever $J$ is a nonnegative continuous convex function. The robust recovery with linear rate of problems (\ref{robu-prob}) and (\ref{regu-prob}) has been discussed in \cite{CRPW,FNT,GHS,ZYY}.

	A necessary and sufficient condition for the uniqueness of solutions to problem $(P_J)$ and a sufficient condition for robust recovery with a linear rate by solving problem (\ref{robu-prob}) via the geometric notion of descent cone have been proposed in \cite{CRPW} for a gauge $J$. Recently, Fadili, Nghia, and Tran \cite{FNT} have proved that for a non-negative continuous convex function, the uniqueness of solutions to problem $(P_J)$ is equivalent to robust recovery for both problems (\ref{robu-prob}) and (\ref{regu-prob}), and that a sharp solution guarantees robust recovery with a linear rate.

	In this paper, we propose a unified approach to investigate solution uniqueness and sharpness in sparse regularization with a closed gauge from a dual perspective. We formulate a Lagrange dual problem $(D_J)$, which is a convex optimization problem with an inequality constraint, for the problem $(P_J)$. However, in general, strong duality cannot hold due to the degeneracy of problem $(P_J)$. It is observed that the Lagrange dual problem of $(D_J)$ is actually the original problem $(P_J)$, and Slater's condition holds, leading to strong duality. Consequently, the solution set for problem $(P_J)$ is the Lagrange multiplier set for problem $(D_J)$, and is same for any solution to $(D_J)$. 
	
	This article is organized as follows. In the next section, we will present some preliminary results, particularly concerning robust recovery with a linear rate. 
	In Section $3$, we first give a characterization of the unique solution via the so-called radial cone of the unit sublevel set. We then establish characterizations of solution uniqueness to problem $(P_J)$ and the calmness properties of two types of solution mappings, and prove the equivalence of these characterizations under the assumption of metric subregularity of the subdifferentials. The $\ell_1/\ell_2$ group Lasso regularizer and the nuclear norm satisfy such conditions. Additionally, we present some sufficient and necessary conditions for a sharp solution, which ensures robust recovery with a linear rate, and also imply that the solution mapping is locally upper Lipschitzian. 
	In Section $4$, we will present some applications for sparse analysis regularization, weighted sorted $\ell_1$-norm sparse regularization, nuclear norm optimization, conic gauge optimization, and semidefinite conic gauge optimization. Finally, a conclusion is given in the last section.

	\section{Preliminaries}
	Throughout the paper, given vectors $x, y \in \R^n$, we denote the inner product by $\langle x, y \rangle,$ the Euclidean norm by $\| x\|$, the $\ell_1$ norm by $\|x \|_1$, the $\ell_{\infty}$ norm by $\|x\|_{\infty}$.  We denote the closed unit ball in Euclidean norm by ${\mathbb B}$, $\ell_{\infty}$ norm by $B_{\infty}$, and  the closed ball centered at $x$ with radius $r>0$ by $B_{r}(x)$. Given a matrix $A \in \R^{m \times n}$, $A^T  $, ${\rm Im}A$, and ${\rm Ker}A$ denote the transpose, the image, and the kernel of $A$, respectively.
	
	The notations and concepts of convex analysis that we employ are standard \cite{HUL,Rockafellar,RW}. Given a nonempty set $S \subset \R^n$, we use ${\rm aff } S$ to represent its affine hull, ${\rm cone} S$ for its conical hull, ${\rm cl } S$ for its closure, ${\rm int} S$ for its interior, and ${\rm ri} S$ for its relative interior. 
	The indicator function of the set $S$, denoted by $I_S$, is defined by $I_S(x) = 0$ if $x \in S$ and $+ \infty$ otherwise. 
	The support function of the set $S$, denoted by $\sigma_S$, is defined by $\sigma_S(\cdot) := \sup_{s\in S}\langle s, \cdot \rangle$.  
	The distance function from a point $x$ to the set $S$ is defined by ${\rm d}(x,S):=\inf\{\|x-y\| \colon y \in S\}$.
	
	Let $f \colon \R^n \to \R \cup \{\pm \infty\}$ be a function. The domain and epigraph of $f$ are defined as
	${\rm dom}f:= \{ x \in \R^n \colon f(x) < \infty\}$ and ${\rm epi}f := \{ (x,r) \in \R^n \times \R \colon f(x) \le r\}.$ $f$ is proper if ${\rm dom}f \ne \emptyset$  and $f(x) > -\infty $ for all $x \in {\rm dom}f$. $f$ is convex and lower semi-continuous (lsc) if ${\rm epi}f$ is convex and closed, respectively. The directional derivative of $f$ at $\bar x$ in the direction $h$ is defined as
	\begin{equation*}
		f^{\prime}(\bar x,h) := \lim_{t \downarrow 0} \frac{f(\bar x + t h) -f(\bar x)}{t},
	\end{equation*}
	and the subderivative of $f$ at $\bar x$ in the direction $h$ is defined as
	\begin{equation*}
		{\rm d} f(\bar x)(h) := \liminf_{ {t \downarrow 0} \atop{ h^{\prime} \to h}} \frac{f(\bar x + t h^{\prime})-f(\bar x)}{t}.
	\end{equation*}
	The subdifferential of $f$ at $\bar x \in {\rm dom} f$ is defined by
	\begin{equation*}
		\partial f(\bar x) : = \{s \in \R^n \colon  f(x) \ge  f(\bar x) + \langle s, x-\bar x \rangle  \mbox{ for all } x \in \R^n\}.
	\end{equation*}
	We denote $f^* \colon \R^n \to \R \cup \{\pm \infty\}$ by the Fenchel conjugate of $f$:
	$$
	f^*(v):=\sup\{\langle v,x \rangle-f(x) \colon x \in \R^n\} \mbox{ for } v \in \R^n.
	$$
	When $C$ is a nonempty closed convex set, the subdifferential of $I_C$ at $\bar x\in C$ is the normal cone to $C$ at $\bar x$, which is defined as
	\begin{equation*}
		N_C(\bar x):= \{s \in \R^n \colon \langle s, x-\bar x \rangle \le 0 \mbox{ for all } x \in C\}.
	\end{equation*}
	Given a convex cone $K \subset \R^n$, the polar of $K$ is defined as the cone
	\begin{equation*}
		K^{o} := \{ s \in \R^n \colon \langle s, x \rangle \le 0 \mbox{ for all } x \in K \}.
	\end{equation*}
	The tangent cone to $C$ at $\bar x \in C$ is defined by $T_C(\bar x) = N_C(\bar x)^{o}$. The cone generated by $C-\bar x$ is 
	$${\rm cone}(C -\bar x) = \bigcup\limits_{\lambda > 0} \lambda (C-\bar x)= \{ h \in \R^n \colon \exists  t^* > 0, \forall t \in [0,t^*], \bar x+t h \in C \} =:{\cal R}_C(\bar x),$$
	which is the radial cone of $C$ at $\bar x$.
	
	Consider a set-valued mapping $F\colon \R^n \rightrightarrows \R^m$ with its domain and graph given by 
	\begin{equation*}
		{\rm dom}F := \{x \in \R^n \colon F(x) \ne \emptyset \}, \mbox{ and } {\rm gph}F:=\{(x,y) \in \R^n \times \R^m \colon y \in F(x)\}.
	\end{equation*} 
	The graphical derivative of $F$ at $(\bar x, \bar y) \in {\rm gph}F$ is $DF(\bar x, \bar y)\colon \R^n \rightrightarrows \R^m$ defined by
	\begin{equation*}
		DF(\bar x, \bar y)(u) = \{ v \in \R^m \colon (u, v) \in T_{{\rm gph}F}(\bar x, \bar y)\}  \;\; \mbox{ for all } u \in \R^n.
	\end{equation*}
	The mapping $F$ is metrically regular around $(\bar x,\bar y) \in {\rm gph}F$ if there exists $\kappa \ge 0$ together with neighborhoods $U$ of $\bar x$ and $V$ of $\bar y$ such that
	\begin{equation}\label{metric-regular}
		{\rm d} (x, F^{-1}(y)) \le \kappa {\rm d} (y, F(x))  \;\; \mbox{ for all } (x,y) \in U \times V. 
	\end{equation} 
	If (\ref{metric-regular}) is satisfied for fixed $y=\bar y$, we say that $F$ is metrically subregular around $(\bar x,\bar y).$  We say that $F$ is strongly metrically subregular at $(\bar x, \bar y)$ if there exists $\kappa \ge 0$ and a neighborhood $U$ of $\bar x$ such that 
	\begin{equation*}
		\|x -\bar x\| \le \kappa {\rm d}(\bar  y,F(x))\;\; \mbox{  for all } x \in U.
	\end{equation*}
	$F$ is calm at $(\bar x, \bar y) \in {\rm gph}F$ with modulus $\kappa \ge 0$ if there exist neighborhoods $U$ of $\bar x$ and $V$ of $\bar y$ such that
	\begin{equation*}
		F(x) \cap V \subset F(\bar x) + \kappa \|x-\bar x\|{\mathbb B}\;\; \mbox{ for all } x \in U,
	\end{equation*}
	where $\mathbb{B}$ denotes the closed unit ball of $\R^m$. Furthermore, we say $F$ is isolated calm at $(\bar x, \bar y) \in {\rm gph}F$ with modulus $\kappa \ge 0$ if $F(\bar x)= \bar y$ and $F$ is calm at $(\bar x, \bar y) \in {\rm gph}F$ with modulus $\kappa$. It is well known (see \cite{DR}) that the calmness and isolated calmness of $F$ at $(\bar x,\bar y)$ are equivalent to the metric subregularity and strong metric subregularity of the inverse mapping $F^{-1}$ at $(\bar y,\bar x)$, respectively.
	
	A gauge is a function $J \colon \R^n \to \R \cup\{+ \infty\}$ that is convex, nonnegative, positively homogeneous of degree one (i.e., $J(\alpha x) = \alpha J(x)$ for all $\alpha > 0$ and $x \in \R^n)$, and $J(0) = 0$ \cite[p.128]{Rockafellar}. The unit sublevel set of a gauge $J$, namely
	\begin{equation*}
		B:= \{x \in \R^n \colon J(x) \le  1\},
	\end{equation*}
	is therefore a convex set that contains the origin. If $J$ is lsc, then it coincides with the
	Minkowski function associated with $B$, i.e., $J = \gamma_B$, where
	\begin{equation*}
		\gamma_{B}(x) = \inf\{\lambda > 0 \colon x \in \lambda B\}.
	\end{equation*}
	The polar of a set $C \subset \R^n$ is the closed convex set defined as
	\begin{equation*}
		C^{o}= \{z \in \R^n \colon \langle z,x \rangle \le 1 \mbox{ for all } x \in C\}.
	\end{equation*}
	The polar of a gauge $J$ is the closed gauge $J^{o}$ defined by
	\begin{equation*}
		J^{o}(z) =\inf \{ \mu \ge 0 \colon \langle z, x\rangle \le \mu J(x) \mbox{ for all } x \in \R^n\}.
	\end{equation*}
	The unit sublevel set of a gauge $J^o$ is $B^o$, and the conjugate of $J$ is $J^* = I_{B^o}.$   
	
	We say $\bar x$ is a sharp minimizer of $f$  if there exist some $\kappa, \delta >0$ such that
	\begin{equation*} \label{sharp-solu}
		f(x) - f(\bar x ) \ge \kappa \|x -\bar x\|  \mbox{ for all } \|x - \bar x  \| \le \delta.
	\end{equation*}
	This concept was independently introduced by Cromme \cite{Cromme} and Polyak \cite{Polyak}. It is a strong condition for ensuring the uniqueness of solutions to a convex function. If the solution is not isolated, meaning the solution set $S$ is not a singleton, a much looser condition is that there exist some $\kappa, \delta >0$ such that
	\begin{equation*}\label{w-grow}
		f(x) - f(\bar x) \ge \kappa {\rm d}(x,S)    \mbox{ for all }  \|x - \bar x  \| \le \delta.
	\end{equation*}
	This notion proposed by Burke and Ferris \cite{BF} is  known as weak sharp minima.
	
	\begin{lem}\label{sharp-sol}(\cite[Theorem 2.2]{BF}) 
		$x^*$ is a sharp  minimizer of a proper lsc convex function $f$ if and only if  $0 \in {\rm int} \partial f(x^*).$
	\end{lem}
	
	It was proved in \cite[Corollary 3.6]{BF} that the solution set of a proper polyhedral convex function (if nonempty) is a set of weak sharp minima. This implies that uniqueness of solutions is equivalent to sharp minima in the case of convex polyhedral problems,  such as problem ($P_J$) with $\ell_1$ or $\ell_{\infty}$ norms.

	When $J$ is a continuous proper convex function or a convex polyhedral function,  a feasible point $x_0$ is a sharp solution of ($P_J$) if and only if $ 0 \in {\rm int} \partial(J+ I_{\mathcal{X}})(x_0)= {\rm int}[\partial J(x_0)+{\rm Im}A^T ]$ by Lemma \ref{sharp-sol}.

	Suppose that $J \colon \R^n \to \bar{\R}$ is a nonnegative lsc function but not necessarily convex. We denote by $J_{\infty} \colon \R^n \to  \bar{\R}$ the asymptotic (or horizon) function  associated with $J$, which is defined by
	\begin{equation*}
		J_{\infty}(w) = \liminf_{w^{\prime} \to w,t \to \infty} \frac{J(tw^{\prime})}{t},
	\end{equation*}
	and denote by ${\rm Ker} J_{\infty} := \{w\in \R^n \colon J_{\infty}(w)=0\}$ the kernel of $J_{\infty}$.
	Throughout this section, we assume that $J$ satisfies
	\begin{equation} \label{coer-cond}
		{\rm Ker} J_{\infty} \cap  {\rm Ker } A = \{0\}.
	\end{equation}
	It is well known that this condition ensures that problem ($P_J$) has a nonempty compact set
	of minimizers and the converse remains true whenever $J$ is convex (see \cite[Corollary 3.1.2]{AT}).
	
	For $J$ a nonnegative lsc function satisfying (\ref{coer-cond}), \cite[Proposition 3.1]{FNT} proves that solution uniqueness is sufficient for robust recovery, and the converse remains true if $J$ is convex and ${\rm dom} J = \R^n$. However, solution uniqueness is not enough to guarantee robust recovery with a linear rate for a general regularizer. The following result provides a sufficient condition for robust recovery with a linear rate by solving the problem (\ref{robu-prob}). Denote by $D_J(x^*):={\rm cone}\{x-x^* \colon J(x) \le J(x^*)\}$ the descent cone of $J$ at $x^*$ and $\tilde{C}_J(x^*) := \{h\in \R^n \colon {\rm d}J(x^*)(h) \le 0 \}$ the critical cone of $J$ at $x^*$.
	
	The following result was proved in \cite[Proposition 2.1, Proposition 2.2]{CRPW} for a gauge and the proof remains true for a proper convex function.  
	\begin{lem} \label{lem2.2} (\cite[Proposition 2.1, Proposition 2.2]{CRPW})  Let $J$ be a  proper lsc convex function. \\
		{\rm (i)} $x_0$ is the unique solution to problem $(P_J)$ if and only if ${\rm Ker A} \cap D_J(x_0) = \{0\}$;\\
		{\rm (ii)} Suppose that $\|Ax_0 - b\| \le \delta$ and there exists some $\alpha > 0$ such that
		\begin{equation}\label{robu-recov-cond} 
			\|Ah\| \ge \alpha \|h\| \mbox{ for all }  h \in D_J(x_0).
		\end{equation}
		Then any solution $x^{\delta}$ to problem (\ref{robu-prob}) satisfies
		\begin{equation*}
			\|x^{\delta} - x_0\| \le \frac{2\delta}{\alpha}.
		\end{equation*}
	\end{lem}

	The supremum of all $\alpha > 0$ in (\ref{robu-recov-cond}) is called the minimum conic singular value of $A$ with respect to $D_J(x_0)$, denoted by
	\begin{equation*}
		\lambda_{\rm min}(A,D_J(x_0)) :=\inf \{\frac{\|Aw\|}{\|w\|}\colon w \in D_J(x_0) \setminus \{0\}\}.
	\end{equation*}
	It was proved in \cite{CRPW,tropp} that $\lambda_{\rm min}(A,D_J(x_0))> 0$ with high probability when the number of measurements $m$ is larger than a function of the Gaussian widths of $D_J(x_0)$ for a standard normal Gaussian matrix $A$.
	
	\begin{lem}\label{lem2.3} (Sharp minima for robust recovery)\cite[Corollary 3.11]{FNT} Suppose that $J$ is a nonnegative lsc proper convex function and locally Lipschitz continuous around $x_0$ with modulus $L$. If $x_0$ is a sharp solution of ($P_J$) with constant $\kappa > 0$, then for any $\delta > 0$, the following statements hold:\\
		{\rm (i)} Any solution $x^{\delta}$ to problem (\ref{robu-prob}) with $\|Ax_0 -b \| \le \delta$ satisfies
		\begin{equation*}
			\|x^{\delta} - x_0\| \le \frac{2(L + \kappa)\|A^{\dagger}\|}{\kappa} \delta,
		\end{equation*}
		where $A^{\dagger}$ denotes the Moore-Penrose generalized inverse of $A$ (see \cite{Meyer}).\\
		{\rm (ii)} For any $c_1 > 0$ and $\mu = c_1\delta$, any minimizer $x^{\mu}$ to problem (\ref{regu-prob}) with $\|Ax_0 -b \| \le \delta$ satisfies
		\begin{equation*}
			\|x^{\mu} - x_0\| \le \frac{c_1}{2\kappa} \Bigg( \frac{1}{c_1} + (\kappa+L)\|A^{\dagger}\| \Bigg)^2 \delta.
		\end{equation*}
	\end{lem}
	
	\medskip
	\noindent {\bf Remark 2.1.}
	When $J$ is convex and continuous,	\cite[Proposition 3.8]{FNT} proves that $x_0$ is a sharp solution of problem $(P_J)$ if and only if
	${\rm Ker} A \cap C_J(x_0) = \{0\}$, where $C_J(x_0) := \{h\in \R^n \colon J^{\prime}(x_0,h) \le 0\}$. It follows from  \cite[Proposition 8.21]{RW} or \cite[Proposition 2.126]{Bonnans-Shapiro} that ${\rm d}J(x_0)(\cdot) = J^{\prime}(x_0,\cdot)$, and then $C_J(x_0) =\tilde{C}_J(x_0)$.  
	Since $C_J(x_0)$ is a closed cone, it is easy to see that  ${\rm Ker}A \cap C_J(x_0) = \{0\}$ if and only if there exists $\alpha > 0$ such that 
	\begin{equation*}\label{sharp-suff-cond}
		\|Ah\| \ge \alpha \|h\| \mbox{ for all }  h \in C_J(x_0).
	\end{equation*}
	This is equivalent to (\ref{robu-recov-cond}) whenever $C_{J}(x_0) = {\rm cl}D_{J}(x_0)$. 
	
	\medskip
	\noindent {\bf Remark 2.2.} Given $x^* \in \R^n$, denoted by $I ={\rm supp}(x^*):= \{ i \in \{ 1,2, \ldots, n\} \colon x^*_i \ne 0\}$ the support set of $x^*$ and $I^c$ the complementary of $I$. $h_I \in \R^n$ satisfies $(h_I)_i =h_i$ for $i \in I$ and zero elsewhere.   For $J(\cdot) = \ell_1(\cdot):=\|\cdot\|_1$, Cand\'{e}s, Romberg, and Tao \cite{CRT} prove that $D_{\ell_1}(x^*) \subset Z_I := \{h \in \R^n \colon \|h_{I^c}\|_1 \le \|h_I\|_1\}$ (it follows clearly that $C_{\ell_1}(x^*) \subset Z_I$ since $C_{\ell_1}(x^*) = {\rm cl}D_{\ell_1}(x^*)$ and $Z_I$ is closed), and 
	\begin{equation*}
		\|Ah\| \ge \alpha \|h\| \mbox{ for all }  h \in Z_I, 
	\end{equation*}
	or equivalently,  ${\rm Ker}A \cap Z_I =\{0\}$ holds provided $\delta_{3k} + 3 \delta_{4k} < 2$, where $\delta_k$ is the restricted isometry constant (RIC) for $A$ and $k = {\rm card}(I)$ is the cardinality of $I$. We also observe that ${\rm Ker}A \cap Z_I =\{0\}$ if and only if $\|h_{I^c}\|_1 > \|h_I\|_1$ for all $h \in {\rm Ker }A \setminus \{0\}.$   \cite[Lemma 2.2]{Candes} shows that the later holds provided $\delta_{2k} < \sqrt{2}-1.$ The above inequality can be equivalently written as $\|h\|_1 > 2\|h_{I}\|_1$ for all $h \in {\rm Ker }A \setminus \{0\},$ or equivalently, $ P_1(I,A)£º=\max\limits_{h \in {\rm Ker}A\setminus \{0\}} \frac{\|h_I\|_1}{ \|h\|_1} < \frac{1}{2},$ which is used in \cite{FR} to guarantee the uniqueness of solution to $(P_1)$.
	
	\section{Uniqueness and Sharpness for Sparse Regularization with a Gauge}
	Unless otherwise specified,  we assume that $J$ is a lsc gauge in this section.  Consider the problem
	$$ \min_{x} J(x) \;\; \mbox{ s.t. } \;\;  Ax=b_0 , \leqno(P_J)$$
	and its Lagrange dual problem is
	$$ \max_{y} \langle y,b_0 \rangle  \;\; \mbox{ s.t. } \;\; J^{o}(A^T  y)\le 1 \Leftrightarrow A^T   y \in B^o. \leqno(D_J)$$
	If the row vectors of $A$ is not linearly independent, the strong duality between problem $(P_J)$ and $(D_J)$ does not necessarily hold. It is easy to show that the Lagrange dual problem to $(D_J)$ is problem $(P_J)$. In fact, the Lagrange function of $(D_J)$ is
	$L(y,x) = \langle y,b_0 \rangle - \langle x, A^T  y\rangle$ and its Lagrange dual problem is given by
	\begin{align*}
		&\min_{x} \{\sup_{y}L(y,x) + \sigma_{B^o}(x) \} \nonumber\\
		=& \min_{x} \{J(x): Ax = b_0\}.
	\end{align*}
	
	If the objective function $J$ is not positively homogeneous, the Lagrange dual problem to $(D_J)$ is not necessarily problem $(P_J)$. For instance,
	consider the problem
	\begin{equation}\label{prim-1}
		\min_x \frac{1}{2}\|x-\bar x\|^2  \;\;\mbox{  s.t. } \;\; Ax  =0,\end{equation}
	and its Lagrange dual problem
	\begin{equation} \label{dual-1}
		\max_{y} -\frac{1}{2}\|A^T   y \|^2 + \langle y, A\bar x\rangle.\end{equation}
	The Lagrange dual problem to problem (\ref{dual-1}) is
	\begin{equation*}\label{dd-1}
		\min_x \frac{1}{2}\|x\|^2  \;\;\mbox{  s.t. } \;\;  Ax  =A\bar x,
	\end{equation*}
	which is not the primal problem (\ref{prim-1}).
	
	Since the weak duality between problem $(P_J)$ and $(D_J)$ always holds, one has ${\rm val}(D_J) \le {\rm val}(P_J)$, and since $y=0$ is feasible for $(D_J)$, one has ${\rm val}(D_J) \ge 0$. So ${\rm val}(D_J)$ must be finite.  Since $(D_J)$ satisfies Slater's condition, by \cite[Theorem 3.4, 3.6]{Bonnans-Shapiro}, the strong duality between problem $(D_J)$ and $(P_J)$ holds.
	
	\begin{lem} \label{Lem 3.1}
		There is no duality gap between problem $(D_J)$ and $(P_J)$, and the Lagrange multiplier set for $(D_J)$ is a nonempty, convex, compact subset and same for any optimal solution of $(D_J)$, and coincides with the optimal solution set of $(P_J)$.
	\end{lem}
	
	Therefore, solution uniqueness to problem $(P_J)$ is equivalent to the uniqueness of Lagrange multipliers for problem $(D_J)$ for any optimal solution of problem $(D_J)$. This is a key basis for our following analysis.
	
	Mangasarian \cite{Manga} provides several sufficient and necessary conditions for the uniqueness of the Lagrange multiplier in linear programming. Kyparisis \cite{Kyparisis} subsequently proves that the uniqueness of the Lagrange multiplier in nonlinear programming is equivalent to the strict Mangasarian-Fromovitz constraint qualification. Shapiro \cite{Shapiro} offers a sufficient and necessary condition for the uniqueness of the Lagrange multiplier in a differentiable optimization problem with conic constraint in Banach spaces (for finite dimensional spaces, see also \cite[Remark 4.5]{MBSE}). Additionally, Bonnans and Shapiro \cite[Theorem 4.47]{Bonnans-Shapiro}  establish that the strict Robinson's condition is sufficient for the uniqueness of the Lagrange multiplier in a differentiable optimization problem with a closed convex set constraint in Banach spaces, and the converse implication remains true under certain conditions where some cones are closed.
	
	\medskip
	\noindent{\bf Assumption 3.1}
	Let the gauge $J$ and the set $\mathcal{X}=\{ x \in \R^n \colon Ax= b_0\}$ satisfy properties ${\rm (A)}$ or ${\rm (B)}$ as follows:\\
	{\rm (A)} $\mathcal{X} \cap {\rm ri(dom}J) \ne \emptyset$.\\
	{\rm (B)} $J$ is polyhedral and $\mathcal{X} \cap {\rm dom}J \ne \emptyset$.
	
	\medskip
	By Lemma \ref{Lem 3.1}, we have the following optimality condition: a point $y_0$ is an optimal solution of $(D_J)$ if and only if there exists a Lagrange multiplier $x$, together with $y_0$ satisfying (see \cite[Theorem 3.4]{Bonnans-Shapiro})
	\begin{equation} \label{optim-condition}
		Ax=b_0 \,\,\mbox{and}\,\, x \in N_{B^o}(A^T  y_0),
	\end{equation}
	or equivalently,
	\begin{equation} \label{optim-condition-1}
		Ax = b_0 \,\,\mbox{and}\,\,  A^T  y_0 \in \partial \sigma_{B^o}(x) = \partial J(x) = B^o \cap \{z \colon \langle x, z \rangle = J(x) \}.
	\end{equation}
	Conversely, if $x_0$ is an optimal solution to problem $(P_J)$, then $0 \in \partial(J+ I_{\mathcal{X}})(x_0)$ and by subdifferential formula \cite[Theorem 23.8]{Rockafellar}, there exists $y_0$ such that $A^T  y_0 \in \partial J(x_0)$ provided Assumption 3.1 holds, which implies that $y_0$ is an optimal solution to $(D_J)$.
	
	Define a set-valued mapping
	\begin{eqnarray*}
		G(x,y):=	  
		\begin{bmatrix}
			-A^T  y\\
			Ax - b_0
		\end{bmatrix}
		+ 
		\begin{bmatrix}
			\partial J(x)\\
			0
		\end{bmatrix}, 		
	\end{eqnarray*}
	and solution mappings
	$S_y(v,w) : = \{x \in \R^{n} \colon  (v,w) \in  G(x,y)\}$
	and 
	$$
	S_{y,0}(w) : = \{x \in \R^{n} \colon  (0,w) \in  G(x,y)\}=\{x \in \partial J^*(A^T  y) \colon w=Ax - b_0\}.
	$$
	
	Clearly, if $S_y(0,0) = S_{y,0}(0)\ne \emptyset $, then $y$ is an optimal solution to problem $(D_J)$, $S_y(0,0)$ is the multiplier set of problem $(D_J)$ and coincides with the solution set of problem $(P_J)$.

	\begin{thm}\label{Thm 3.1} Suppose that Assumption 3.1 holds. A feasible point $x_0$ with $J(x_0)$ finite is the unique solution to  problem $(P_J)$ if and only if there exists $y_0$ satisfying $A^T y_0 \in \partial J(x_0)$ such that
		\begin{equation} \label{cond-unique-solution}
			{\rm Ker}A \cap (N_{B^o}(A^T y_0)+ [x_0]) = \{0\},
		\end{equation}
		where $[x_0]$ denotes the one dimensional space generated by $x_0$. 
	\end{thm}
	{\it Proof}  Under each of assumptions, if $x_0$ is an optimal solution to problem $(P_J)$, then there exists $y_0$ satisfying  $A^T  y_0 \in \partial J(x_0)$, which is equivalent to $y_0$ being an optimal solution to $(D_J)$.  A vector $x$ is an optimal solution to problem $(P_J)$ if and only if it is a Lagrange multiplier of $(D_J)$ for any solution to $(D_J)$. This means that  $(x,y_0)$ satisfies (\ref{optim-condition}), which is equivalent to
	\begin{equation*}
		A(x-x_0)=0 \mbox{ and }  x-x_0 \in N_{B^o}(A^T  y_0)-x_0.
	\end{equation*}
	It follows that $x_0$ is the unique solution to problem $(P_J)$ if and only if
	\begin{equation*}
		{\rm Ker}A \cap {\rm cone}(N_{B^o}(A^T  y_0)-x_0) = \{0\}.
	\end{equation*}
	This is exactly (\ref{cond-unique-solution}). 
	\qed
	
	\medskip 
	\noindent {\bf Remark 3.1.}
	When $J^o (A^T   y_0) = 1$,  it follows from \cite[Lemma 26.17]{BC} that $x_0$ is the unique solution to problem $(P_J)$ if and only if
	\begin{equation*} 
		{\rm Ker}A \cap [(N_{{\rm dom}J^o}(A^T  y_0) \cup {\rm cone}\{x \in B \colon \langle A^T   y_0, x \rangle =1 \})+ [x_0]] = \{0\}.   
	\end{equation*}

	Note that $J^*= I_{B^o}$, then  ${\rm cone}(N_{B^o}(A^T  y_0)-x_0)=  {\cal R}_{N_{B^o}(A^T  y_0)}(x_0)= {\cal R}_{\partial J^*(A^T  y_0)}(x_0).$ Therefore, condition (\ref{cond-unique-solution}) is equivalent to 
	\begin{equation*}{\rm Ker}A \cap {\cal R}_{\partial J^*(A^T  y_0)}(x_0) = \{0\}.
	\end{equation*} 
	The conclusion of Theorem 3.1, for the case when $J = \|\cdot\|_*$ is the nuclear norm, has already been established in a previous work \cite[Lemma 3.2]{LPB}. Recently, \cite[Theorem 4.5]{FNP-1} has proved the same conclusion for a continuous convex function $J$ by using the equality $ {\rm Ker}A \cap D_J(x_0) ={\rm Ker}A \cap {\cal R}_{\partial J^*(A^T  y_0)}(x_0)$ (as presented in \cite[Proposition 4.4]{FNP-1}).
	
	Building upon Theorem \ref{Thm 3.1}, we are able to derive the following result.

	\begin{thm}\label{Thm 3.2} Suppose that Assumption 3.1 holds. Let $x_0 \in \R^n$ be such that there exists $y_0$ satisfying $A^T  y_0\in \partial J(x_0).$  Then the following conditions are equivalent:\\
		{\rm (i)}
		$x_0$ is  the unique solution to problem $(P_J)$ and $S_{y_0,0}(\cdot)$ is calm at $(0,x_0)$;\\
		{\rm (ii)}
		$x_0$ is  the unique solution to problem $(P_J)$ and there exists $\kappa >0$ and $\epsilon >0$ such that
		\begin{equation*}
			{\rm d}(x, S_{y_0,0}(0)) \le \kappa  \| Ax -b_0\|   \mbox{ for all } x \in  \partial J^*(A^T  y_0) \cap B_{\epsilon}(x_0); 
		\end{equation*}
		{\rm (iii)}
		\begin{equation*} \label{cond-unique-solution-2-1}
			{\rm Ker}A \cap T_{\partial J^*(A^T  y_0)}(x_0) = \{0\};
		\end{equation*}
		{\rm (iv)}
		\begin{equation*} \label{cond-unique-solution-1}
			{\rm Ker}A \cap {\rm cl}(N_{B^o}(A^T  y_0)+[x_0]) = \{0\};
		\end{equation*}
		{\rm (v)}
		\begin{equation*} \label{cond-unique-solution-2}
			{\rm Im}A^T   + T_{B^o}(A^T  y_0) \cap [x_0]^{\bot} = \R^n.
		\end{equation*}
		
		Moreover,	if $N_{B^{o}}(A^T  y_0)+ [x_0]$ is closed, especially in the case where $B$ is polyhedral, then the above conditions are equivalent to $x_0$ being the unique solution to problem $(P_J)$.
	\end{thm}

	{\it Proof}  It follows from the definition of calmness of $S_{y,0}$ at $(0,x_0)$ that  implication ${\rm (i)} \Leftrightarrow {\rm (ii)}$ holds.
	Set $L(x) := (Ax-b_0,x)$ for any $x \in \R^n$. Then we have 
	${\rm gph}S_{y_0,0} = L(\partial J^*(A^T  y_0))$.
	By \cite[Proposition A.5.3.1]{HUL}, we obtain that 
	\begin{eqnarray*}
		T_{{\rm gph}S_{y_0,0}}(0,x_0) &=& {\rm cl} \{ (\xi, \eta) \colon \xi = A\eta, \eta \in T_{\partial J^*(A^T  y_0)} (x_0) \}\\
		&=& \{ (\xi, \eta) \colon \xi = A\eta, \eta \in T_{\partial J^*(A^T  y_0)} (x_0) \}.
	\end{eqnarray*}  
	It follows from the definition of graphical derivative that
	$$ DS_{y_0,0}(0,x_0)(\xi) = A^{-1}(\xi) \cap T_{\partial J^*(A^T  y_0)} (x_0).$$
	By \cite[Proposition 4.1]{Levy1996}, $S_{y_0,0}(\cdot)$ is isolated calm at $(0,x_0)$, i.e., $S_{y_0,0}(0)=\{x_0\}$ and $S_{y_0,0}(\cdot)$ is calm at $(0,x_0)$, if and only if  $DS_{y_0,0}(0,x_0)(0)=\{0\}$, which is exactly condition {\rm (iii)}. This implies that  condition  ${\rm (i)} \Leftrightarrow {\rm (iii)}$.

	Since $N_{B^o}(A^T  y_0)+[x_0] = {\cal R}_{\partial J^*(A^T  y_0)}(x_0)$ and $T_{\partial J^*(A^T  y_0)}(x_0)={\rm cl}{\cal R}_{\partial J^*(A^T  y_0)}(x_0)$, implication ${\rm (iii)} \Leftrightarrow {\rm (iv)}$ holds.
	Taking the polar cone of both sides, we can establish that condition {\rm (iv)} is equivalent to  
	\begin{equation} \label{strong-robinson} 
		{\rm cl}[{\rm Im}A^T   + T_{B^o}(A^T  y_0) \cap [x_0]^{\bot}] = \R^n.
	\end{equation}
	Since for any convex set $C$,  ${\rm ri}({\rm cl}C) ={\rm ri}C \subset C$. Therefore, we can deduce the equivalence between (\ref{strong-robinson}) and condition {\rm (v)}.
	
	The implication from condition {\rm (iv)} to (\ref{cond-unique-solution}) is straightforward and the converse implication clearly holds true if $N_{B^{o}}(A^T  y_0)+ [x_0]$ is closed. In the case when $B$ is polyhedral, $T_{B^o}(A^T  y_0)$ is also a polyhedral cone, and then $N_{B^o}(A^T  y_0)+ [x_0]$ is closed due to  Farkas Lemma \cite[Theorem 2.201]{Bonnans-Shapiro}. 
	\qed

	\medskip
	\noindent {\bf Remark 3.2.}  $N_{B^o}(A^T  y_0)+ [x_0] = {\cal R}_{\partial J^*(A^T  y_0)}(x_0)$ is closed whenever $\partial J^*(A^T  y_0)$ is a polyhedral. This condition is automatic when $J$ is a convex piecewise linear-quadratic function \cite[Definition 10.20]{RW}, i.e., ${\rm dom}J$ can be represented as the union of finitely many polyhedral sets, relative to each of which $J(\cdot)$ is given by a convex quadratic expression. Another important example that is not convex piecewise linear-quadratic is the $\ell_1/\ell_2$ group Lasso regularizer defined by
	$$ J(x)=\|x\|_{1,2} := \sum_{S \in \cal J}\|x_S\| \mbox{ for } x \in \R^n,$$
	where $\cal J$ is a partition of $\{1,2 \ldots, n\}$ with $p$ distinct groups (see \cite[Proposition 7]{ZS}). It is also clear that $ {\cal R}_{\partial J^*(A^T  y_0)}(x_0)$ is closed if $x_0 \in {\rm ri} \partial J^*(A^T  y_0)$. Indeed, this condition means that $ {\cal R}_{\partial J^*(A^T  y_0)}(x_0)$ is a subspace and hence is closed.

	The following result can be deduced from \cite[Corollary 4.2]{GO}, for completeness we give the proof. 
	\begin{lem}\label{lemma3.2}
		Let $f \colon \R^n \to (-\infty, + \infty]$ be a proper lsc convex function. Then 
		\begin{equation*}
			T_{\partial f^*(x_0^*)}(x_0) \subset D \partial f^*(x_0^*,x_0)(0), 
		\end{equation*} 
		and the inverse inclusion remains true if  $\partial f$ is metrically subregular around $(x_0,x_0^*).$
	\end{lem}
	{\it Proof}  Take any $h \in T_{\partial f^*(x_0^*)}(x_0)$, then there exist sequences $t_k \downarrow 0$ and $h_k \to h$ such that 
	$x_0 + t_k h_k \in \partial f^*(x_0^*) =\partial f^*(x_0^* + t_k 0).$
	It follows from the definition of graphical derivative that   $h\in	D \partial f^*(x_0^*,x_0)(0).$ 
	
	Conversely, suppose $h\in	D \partial f^*(x_0^*,x_0)(0).$ Then there exist $t_k \downarrow 0,$ $h_k \to h$ and $h_k^* \to 0$ such that 
	$x_0 + t_k h_k \in  \partial f^*(x_0^* + t_k h_k^*).$
	Since  $\partial f$ is metrically subregular at $(x_0, x_0^*)$ if and only if $\partial f^*$ is calm at $(x_0^*,x_0)$, i.e., there exists $\kappa >0$ and neighborhoods $U$ of $x_0^*$ and $V$ of $x_0$ such that 
	\begin{equation*}
		\partial f^*(x^*) \cap V \subset \partial f^*(x_0^*) + \kappa \|x^* - x_0^*\|{\mathbb B}  \mbox{ for all }  x^* \in U.
	\end{equation*} 
	For $k$ sufficiently large, there exists $h_k^{\prime}$ with $\|h_k^{\prime}\|\le 1$ such that 
	$$x_0 + t_k (h_k - \kappa \|h_k^*\|h_k^{\prime})= x_0 + t_k h_k - \kappa \|x_0^*+t_k h_k^*-x_0^*\|h_k^{\prime}\in \partial f^*(x_0^*).$$ 
	Note that $h_k - \kappa \|x_k^*\|h_k^{\prime} \to h$ as $k \to \infty$, by the definition of the contingent cone, it follows that $h \in T_{\partial f^*(x_0^*)}(x_0)$.
	\qed
	
	\medskip 
	It has been proved in \cite{ZS} that the subdifferential of the $\ell_1/\ell_2$ group Lasso regularizer and the nuclear norm are metrically subregular. The following theorem illustrates that the equivalent conditions discussed in Theorem \ref{Thm 3.2} are equivalent to the uniqueness of the solution and the calmness property of the solution mapping $S_{y_0}$ when $\partial J$ is metrically subregular.

	\begin{thm} \label{Thm 3.3} Let $x_0 \in S_{y_0}(0,0)$. Then the following assertions are equivalent:\\
		{\rm (i)} $x_0$ is the unique solution to problem $(P_J)$ and there exists $\kappa >0$ and $\epsilon >0$ such that
		\begin{equation*} \label{metric-subreg}
			{\rm d}(x, S_{y_0}(0,0)) \le \kappa [ \| Ax -b_0\| + {\rm d}(A^T  y_0, \partial J(x))]\mbox{ for all } x \in B_{\epsilon}(x_0); 
		\end{equation*} 
		{\rm (ii)} $x_0$ is the unique solution to problem $(P_J)$  and the mapping $S_{y_0}$ is calm at $((0,0),x_0)$;\\
		{\rm (iii)} 
		\begin{equation*} \label{unique-condi-10}
			{\rm ker} A \cap D \partial J^*(A^T   y_0,x_0)(0) = \{0\}.
		\end{equation*} 
		
		Moreover, if $\partial J$ is metrically subregular at $(x_0,A^T  y_0)$, then the above conditions are equivalent to that in Theorem \ref{Thm 3.2}
	\end{thm}
	
	{\it Proof}  Set $G_{y}(x)= G(x,y)$. It is clear that  $S_{y_0}(v,w) = G_{y_0}^{-1}(v,w)$. Therefore, $S_{y_0}$ is calm at $((0,0),x_0)$ if and only if $G_{y_0}$ is metrically subregular at $(x_0,(0,0)),$ i.e., there exists $\epsilon >0$ and $\kappa>0$ such that 
	\begin{equation*}
		{\rm d}(x, G_{y_0}^{-1}(0,0)) \le \kappa [{\rm d} ((0,0), G_{y_0}(x) )]\mbox{ for all } x \in B_{\epsilon}(x_0).
	\end{equation*}  
	This is exactly condition {\rm (iv)}.  Note that $  \partial J^*(A^T  y_0)= N_{B^o}(A^T  y_0)$ and  $\partial J(x_0)= N_{B^o}^{-1}(x_0)$, according to \cite[Theorem 4.1]{MBSE}, it can be deduced that conditions {\rm (i)} and {\rm (iii)}  are equivalent. By using Lemma \ref{lemma3.2}, it can be concluded that   condition {\rm (iv)}  implies condition {\rm (ii)} from Theorem \ref{Thm 3.2} and the converse implication holds true if $\partial J$ is metrically subregular at $(x_0,A^T  y_0)$.
	\qed

	\medskip  
	\noindent{\bf Remark 3.3.} By using the calculations of  the subdifferential $\partial J^*$ and the tangent cone $T_{{\rm gph} \partial J^*}$,    \cite[Theorem 3.2]{LPB} proves that conditions {\rm (vi)} and {\rm (v)} from Theorem \ref{Thm 3.2} and condition {\rm (iii)} from Theorem \ref{Thm 3.3} are equivalent and sufficient for the solution uniqueness in the case where $J$ is the nuclear norm.
	
	It is easy to show that  if either $\mathcal{X} \cap {\rm ri} \partial J^*(A^Ty_0) \ne \emptyset$ or  $\partial J^*(A^T  y_0)$ is polyhedral, then $\mathcal{X}$ and $\partial J^*(A^T y_0)$ are boundedly linearly regular (see \cite[Fact 3]{ZS}), which implies that $S_{y_0,0}$ is calm at $(0,x_0)$.
	If, in addition, $\partial J$ is metrically subregular at $(x_0,A^T  y_0)$, $S_{y_0}$ is calm at $((0,0),x_0)$ by \cite[Theorem 2]{ZS}.
	Clearly, the $\ell_1/\ell_2$ group Lasso regularizer satisfies both conditions.
	
	\cite[Theorem 3.5]{FNP} proves that condition {\rm (iv)} from Theorem \ref{Thm 3.2}  implies that $x_0$ is a strong solution to $(P_J),$ i.e., there exists $c >0$ and $\delta >0$ such that 
	$$J(x) - J(x_0) \ge c \|x -x_0\|^2  \; \mbox{ for all } x\; \mbox{ with } \|x - x_0\| < \delta \mbox{ and } A x =b_0,$$ 
	whenever $J$ is second order regular at $x_0$ in the sense of \cite[Definition 3.85]{Bonnans-Shapiro} and $\partial J$ is metrically subregular at $x_0$. In particular,  the $\ell_1/\ell_2$ norm and the nuclear norm satisfy the required conditions (see \cite{ZS}). 
	
	The following result gives the characterizations of a sharp solution.
	
	\begin{thm}\label{Thm 3.5}  Suppose that Assumption 3.1 holds.  A feasible point $x_0$ with $J(x_0)$ finite is a sharp solution of problem $(P_J)$ if and only if there exists $y_0$ satisfying $A^T  y_0 \in \partial J(x_0)$ such that any of the equivalent conditions hold:\\
		{\rm (i)}
		\begin{equation*} \label{cond-unique-solution-4}
			{\rm Ker}A \cap  N_{\partial J(x_0)}(A^T  y_0) = \{0\};
		\end{equation*}
		{\rm (ii)}
		\begin{equation*} \label{cond-unique-solution-4-1}
			{\rm Im}A^T   - T_{\partial J(x_0)}(A^T  y_0)=\R^n;
		\end{equation*}
		{\rm (iii)}
		\begin{equation*} \label{cond-unique-solution-5}
			{\rm Im}A^T   - {\rm cone}(\partial J(x_0) -A^T  y_0)=\R^n;
		\end{equation*}
		{\rm (iv)} The strict Robinson's condition for problem $(D_J)$ with respect to $x_0$ holds, i.e.,
		\begin{equation*} \label{cond-unique-solution-6}
			0 \in {\rm int} [{\rm Im}A^T   - \partial J(x_0)];
		\end{equation*}
		{\rm (v)}
		\begin{equation*}\label{cond-unique-solution-7}
			{\rm Im}A^T   - {\rm aff}(\partial J(x_0)) = \R^n \mbox{ and } {\rm Im}A^T   \cap {\rm ri}(\partial J(x_0)) \ne \emptyset;
		\end{equation*}
		{\rm (vi)}
		\begin{equation*} \label{cond-unique-solution-8}
			{\rm Ker}A \cap  ({\rm cone} \partial J(x_0))^o = \{0\};
		\end{equation*}
		{\rm (vii)}
		\begin{equation*} \label{cond-unique-solution-9-1}
			{\rm Ker}A \cap  \tilde{C}_J(x_0) =\{0\};
		\end{equation*}
		{\rm (viii)}
		\begin{equation*} \label{cond-unique-solution-3}
			{\rm Im}A^T   - {\rm cone}(B^{o} -A^T  y_0)  \cap [x_0]^{\bot}=\R^n.
		\end{equation*}
		
		Moreover, the above equivalent conditions imply all equivalent conditions from Theorem \ref{Thm 3.2} and both classes of conditions are equivalent if ${\rm cone }(B^{o} -A^T  y_0)$ is closed. 
	\end{thm}
	{\it Proof}  Since the set
	$${\rm Im}A^T  - {\cal R}_{\partial J(x_0)}(A^T  y_0) ={\rm Im}A^T   - {\rm cone}(\partial J(x_0) -A^T  y_0)$$ 
	has a nonempty relative interior, by \cite[Proposition 2.97]{Bonnans-Shapiro} we have conditions 	{\rm (i)}-{\rm (iii)} are equivalent to each other and to 
	\begin{equation*} \label{str-robinson} 
		0 \in {\rm int} [A^T  y_0 + {\rm Im}A^T   - \partial J(x_0)],
	\end{equation*}
	which is the same as condition {\rm (iv)}. Furthermore, it can be readily observed that condition {\rm (iv)} is equivalent to condition {\rm (v)}.
	
	According to \cite[Proposition 2.16]{Bonnans-Shapiro} (see also  \cite[Proposition 8.30]{RW}), we can obtain 
	\begin{equation*}
		{\rm d}J(x_0)(h)  = \sigma_{\partial J(x_0)}(h).
	\end{equation*}
	Therefore, condition {\rm (vi)} can be regarded as equivalent to condition {\rm (vii)}.
	
	In accordance with formula (2.32) in reference \cite{Bonnans-Shapiro}, it can be deduced that by computing the polar cone of both sides of condition {\rm (vi)},  condition {\rm (vi)} is equivalent to
	\begin{align}\label{cond-unique-solution-8-eq}
		{\rm cl}[{\rm Im}A^T   - {\rm cone}\partial J(x_0)] = {\rm cl}[{\rm Im}A^T   - {\rm clcone}\partial J(x_0)] =\R^n.
	\end{align}
	Since ${\rm Im}A^T   - {\rm cone}\partial J(x_0)$ has a nonempty relative interior,  by \cite[Theorem 2.17]{Bonnans-Shapiro} we have that (\ref{cond-unique-solution-8-eq})  is equivalent to
	\begin{align} \label{cond-unique-solution-9-1-1}
		0 \in {\rm int} [{\rm Im}A^T   - {\rm cone}\partial J(x_0)].
	\end{align}
	It is obvious that condition {\rm (iv)} implies
	(\ref{cond-unique-solution-9-1-1}) and therefore leads to condition {\rm (vi)}. 
	
	It is clear that condition {\rm (vii)} implies that
	\begin{equation*}
		{\rm d}(J + I_{\mathcal{X}})(x_0)(h) \ge  {\rm d}J(x_0)(h) + I_{{\rm ker}A}(h) > 0 \mbox{ for all } h \ne 0,
	\end{equation*}
	which is equivalent to $x_0$ being a sharp solution of problem $(P_J)$ by \cite[Proposition 3.99]{Bonnans-Shapiro}.
	Then under Assumption 3.1, it follows from Lemma \ref{sharp-sol} that $x_0$ is a sharp solution of problem $(P_J)$ if and only if
	$
	0 \in {\rm int} \partial (J + I_{\mathcal{X}})(x_0) = {\rm int}[\partial J(x_0) + {\rm Im}A^T  ],
	$
	precisely as expressed in condition {\rm (vi)}.  
	Therefore, all the equivalent conditions {\rm (i)} - {\rm (vii)}  are sufficient and necessary conditions for $x_0$ being a sharp solution of problem $(P_J)$.
	
	Since $J$ is a gauge, it can be deduced from  (\ref{optim-condition-1}) that $\partial J(x_0) -A^T  y_0 = (B^o - A^T  y_{0}) \cap [x_0]^{\bot}$, and that
	\begin{equation*} \label{radi-cone}
		{\rm cone}(\partial J(x_0)- A^T  y_0) = {\rm cone }(B^{o} -A^T  y_0)  \cap [x_0]^{\bot}.
	\end{equation*}
	Therefore, the equivalence between condition {\rm (viii)} and  condition {\rm (iii)} is clear. 
	
	Clearly, condition {\rm (viii)} implies condition {\rm (v)} from Theorem \ref{Thm 3.2}. 
	When ${\rm cone }(B^{o} -A^T  y_0)$ is closed,   we have 
	$T_{B^o}(A^T  y_0) ={\rm cone }(B^{o} -A^T  y_0)$ and that  {\rm (viii)} is equivalent to
	condition {\rm (v)} from Theorem \ref{Thm 3.2}.
	\qed

	\medskip
	\noindent {\bf Remark 3.4.} 
	From the proof, it is easy to see that the equivalence among conditions {\rm (i)} - {\rm (vii)} remains true when $J$ is a lsc proper convex function.
	
	The condition 	${\rm Im}A^T   - {\rm aff}(\partial J(x_0)) = \R^n$ is equivalent to $A$ being restricted injective on the space $T := [{\rm aff}(\partial J(x_0)- \partial J(x_0))]^{\bot} $.   
	It is proven in \cite[Lemma 3.2, Proposition 5.3]{Gilbert} that condition {\rm (v)} from Theorem \ref{Thm 3.5} is a sufficient and necessary condition for solution uniqueness to problem $(P_J)$ with $J$ as a polyhedral gauge. 
	Based on the subdifferential decomposability, \cite[Corollary 1]{VGFP} proves that condition {\rm (v)} from Theorem \ref{Thm 3.5} is a sufficient condition for $x_0$ being the unique solution to problem $(P_J)$  and \cite[Theorem 4.6]{FNT}  further proves that condition {\rm (v)} from Theorem \ref{Thm 3.5} is a characterization of $x_0$ being a sharp solution of problem  $(P_J)$ in the case of $J$  being a finite-valued nonnegative continuous convex function. 
	When $J$ is the nuclear norm, \cite[Theorem 3.1]{LPB} shows that condition {\rm (vii)} from Theorem \ref{Thm 3.5} is a characterization of the sharp solution and \cite[Lemma 3.2]{LPB} illustrates that the condition {\rm (iv)} from Theorem \ref{Thm 3.2} is not equivalent to condition {\rm (vii)} from Theorem \ref{Thm 3.5}.
	
	Theorem \ref{Thm 3.5} contains several other equivalent conditions. The method is based on the dual approach, which is different from the literature mentioned above, and the function $J$ is not restricted to be polyhedral or continuous.
	
	The following theorem shows that conditions from Theorem \ref{Thm 3.5} imply that the solution mapping $S_{y_0}$ is locally upper Lipschitzian. The proof is in the spirit of \cite [Proposition 4.47]{Bonnans-Shapiro}.
	
	\begin{thm}\label{Prop 3.1} Suppose that Assumption 3.1 holds and that $x_0 \in \R^n$ is a feasible point of problem $(P_J)$. If there exists $y_0$ satisfying $A^T  y_0 \in \partial J(x_0)$ such that any of the equivalent conditions from Theorem \ref{Thm 3.5} are satisfied, then the set-valued mapping $S_{y_0}$ is upper Lipschitzian at $(0,0)$, i.e., there exists $\kappa>0$ such that
		\begin{equation*} \label{upper-lip}
			S_{y_0}(v,w) \subset \{x_0\} + \kappa (\|v\| +\|w\|) {\mathbb B_{\R^n}}, 
		\end{equation*}
		for all $(v,w)$ near $(0,0)$.
	\end{thm}
	
	{\it Proof}  Consider $x \in S_{y_0}(v,w)$, we have $v+A^T  y_0 \in \partial J(x)$ and $w=Ax-b_0$.  Let $x^* \in {\mathbb B_{\R^n}}$ be such that $\frac{1}{2} \|x-x_0\| \le - \langle x-x_0, x^* \rangle$. 
	Suppose that {\rm (iv)} from Theorem \ref{Thm 3.5} is satisfied, then by the Generalized Open Mapping Theorem \cite[Theorem 2.70]{Bonnans-Shapiro}, there exists $\epsilon >0$ such that
	\begin{equation*}
		\epsilon {\mathbb B_{\R^n}} \subset A^T   {\mathbb B_{\R^m}} - (\partial J(x_0) - A^T  y_0).
	\end{equation*}
	Therefore, there exists $\overline{y} \in {\mathbb B_{\R^m}}$ and $s \in \partial J(x_0)$ such that
	$\epsilon x^* = A^T  \overline{y} + A^T  y_0 - s.$  Since  $x \in \partial J^{*} (v+A^T  y_0) = N_{B^o}(v+A^T  y_0)$, we have
	$\langle x, s - v -A^T  y_0 \rangle \le 0$.   This together with  $\partial J(x_0) -A^T  y_0 = (B^o - A^T  y_{0}) \cap [x_0]^{\bot}$  implies that 
	\begin{equation*}
		\langle x-x_0, s -A^T  y_0 \rangle = \langle x, s -A^T  y_0 \rangle\le \langle x, v \rangle \le (\|x-x_0\| + \|x_0\|) \|v\|.
	\end{equation*} 	
	Combining this with $w=Ax-b_0$, we obtain 
	\begin{eqnarray*}
		\frac{1}{2} \epsilon \|x-x_0\| &\le&  - \langle  x-x_0,  \epsilon x^*\rangle = \langle x-x_0, s -A^T  y_0  \rangle - \langle Ax-Ax_0, \overline{y}\rangle \\
		&\le &  (\|x-x_0\| + \|x_0\|) \|v\| + \|w\|.
	\end{eqnarray*}
	When $\|v\| \le \frac{\epsilon}{4}$, it follows that  
	\begin{eqnarray*}
		\|x-x_0\| \le  \frac{4}{\epsilon}(\|x_0\| \|v\|+\|w\| ),
	\end{eqnarray*}
	which completes the proof.
	\qed

	\bigskip
	If the strict Robinson's condition for problem $(D_J)$ with respect to $x_0$ holds, i.e., $0 \in {\rm int}[ {\rm Im}A^T   - \partial J(x_0)]$, then there exists some $\alpha > 0$ such that $B(0,\alpha)  \subset \partial J(x_0) + {\rm Im}A^T  $.  Define
	\begin{equation} \label{sharp-constant}
		\kappa := \sup \{ \alpha >0 \colon B(0,\alpha)  \subset \partial J(x_0) + {\rm Im}A^T  \}.
	\end{equation}
	Then $B(0,\kappa) \subset  \partial J(x_0) + {\rm Im}A^T  $ provided $\partial J(x_0) + {\rm Im}A^T  $ is closed, which is true under Assumption 3.1. Therefore, in this case, $x_0$ is a sharp minimizer of problem $(P_J)$ with modulus $\kappa$ due to \cite[Theorem 2.6]{BF}. As a direct consequence of Lemma \ref{lem2.3} and Theorem \ref{Thm 3.5}, we have the following result.

	\begin{thm}  \label{Thm 3.6} 	Suppose that $J$ is locally Lipschitz continuous with modulus $L$ around $x_0$, which is a feasible point of problem $(P_J)$ . If there exists $y_0$ satisfying $A^T  y_0 \in \partial J(x_0)$ such that any of the equivalent conditions in Theorem \ref{Thm 3.5} holds, then for any $\delta > 0$, the following statements hold:\\
		{\rm (i)} Any solution $x^{\delta}$ to problem (\ref{robu-prob}) with $\|Ax_0 -b \| \le \delta$ satisfies
		\begin{equation*}
			\|x^{\delta} - x_0\| \le \frac{2(L + \kappa)\|A^{\dagger}\|}{\kappa} \delta;
		\end{equation*}\\
		{\rm (ii)} For any $c_1 > 0$ and $\mu = c_1\delta$, any minimizer $x^{\mu}$ to problem (\ref{regu-prob}) with $\|Ax_0 -b \| \le \delta$ satisfies
		\begin{equation*}
			\|x^{\mu} - x_0\| \le \frac{c_1}{2\kappa} \Bigg( \frac{1}{c_1} + (\kappa+L)\|A^{\dagger}\| \Bigg)^2 \delta,
		\end{equation*}
		where $\kappa$ is the constant defined in (\ref{sharp-constant}).
	\end{thm}

	\begin{prop} \label{Prop3.2} Suppose that $J$ is finite at $x^*.$  Then
		\begin{equation*}
			\tilde{C}_J(x^*) = T_{L_{x^*}(J)}(x^*) = {\rm cl}D_J(x^*),
		\end{equation*}
		where $L_{x^*}(J) := \{ x \in \R^n \colon J(x) \le J(x^*)\}$ denotes the lower level set of $J$ for $J(x^*)$.
	\end{prop}
	{\it Proof}  Since $J$ is lsc convex,  it follows from \cite[Proposition 2.55]{Bonnans-Shapiro} that the second equality holds. Additionally, according to \cite[Proposition 2.61]{Bonnans-Shapiro}, we have	$T_{L_{x^*}(J)}(x^*) \subset \tilde{C}_J(x^*)$ and the converse inclusion holds as long as $0 \not\in \partial J(x^*).$

	We next prove that the first equality remains true if $0 \in \partial J(x^*).$ Suppose that  $0 \in \partial J(x^*).$ Then $J(x^*) = \inf_{x}J(x).$ Since $J$ is nonnegative, positively homogeneous, and $J(0) = 0$, one has that $J(x^*)=0.$ It follows that $L_{x^*}(J) = {\rm Ker}J$ is a closed subspace and that
	$T_{L_{x^*}(J)}(x^*) = {\rm Ker}J.$
	Therefore, we obtain 
	\begin{align*} 
		\tilde{C}_J(x^*) &=({\rm cone}\partial J(x^*))^{o}
		= [\partial J(0) \cap \{z \in \R^n \colon \langle z, x^* \rangle = J(x^*)\}]^{o}
		= {\rm cl} [(\partial J(0))^{o} + [x^*]]\\
		&= {\rm cl} [\{ x \in \R^n \colon J(x)=\sigma_{\partial J(0)}(x) \le 0\} +  [x^*]]
		= {\rm cl} [{\rm Ker}J + [x^*]] \subset {\rm Ker}J,
	\end{align*}
	and the proof is complete.
	\qed

	When $J$  is a polyhedral function, it follows from the closeness of $D_J(x^*)$ that $D_J(x^*)=\tilde{C}_J(x^*)$.
	Proposition \ref{Prop3.2} generalizes \cite[Proposition 3.7]{FNT}, where $J$ is the support function of a nonempty compact convex set containing $0$ as a relative interior, to a gauge case.
	
	As a direct consequence of Theorem \ref{Thm 3.5} and Proposition \ref{Prop3.2}, we observe that the condition (\ref{robu-recov-cond}) is equivalent to $x_0$ being a sharp solution of problem $(P_J)$ under the assumption that Assumption 3.1 holds. This leads to a robust recovery with a linear rate due to Lemma \ref{lem2.2}.
	
	\begin{thm}\label{Thm3.6}
		Suppose that Assumption 3.1 holds and that $J$ is a lsc proper convex function.  A feasible point $x_0$ with $J(x_0)$ finite is a sharp solution of problem $(P_J)$ if and only if (\ref{robu-recov-cond}) holds, i.e., there exists some $\alpha >0$ such that
		\begin{equation*}
			\|Ah\| \ge \alpha \|h\| \mbox{ for all } h \in D_J(x_0).
		\end{equation*}
		Moreover, in this case, we have\\
		{\rm (i)} Any solution $x^\delta$ to problem (\ref{robu-prob}) with $\|Ax_0 -b \| \le \delta$, there exists $\alpha > 0$ such that
		\begin{equation} \label{rob-1}
			\|x^{\delta}- x_0 \| \le \frac{2 \delta}{\alpha}.
		\end{equation}
		{\rm (ii)} For any $c_1 > 0$ and $\mu = c_1 \delta$, any solution $x^\mu$ to problem (\ref{regu-prob}) with $\|Ax_0 -b \| \le \delta$, there exists $\alpha > 0$ such that
		\begin{equation*} \label{rob-2}
			\|x^{\mu}- x_0 \| \le \frac{2(1+2c_1 \|y_0\|)}{\alpha}\delta.
		\end{equation*}
	\end{thm}
	{\it Proof}	 According to Theorem \ref{Thm 3.5}, a feasible point $x_0$ with $\partial J(x_0)\ne \emptyset$  is a sharp solution of problem $(P_J)$ if and only if
	\begin{equation*}
		{\rm Ker}A \cap  \tilde{C}_J(x_0) = \{0\}.
	\end{equation*}
	This is due to the fact that $\tilde{C}_J(x_0)$ is a closed cone, which is further equivalent to 
	\begin{equation*}
		\|Ah\| \ge \alpha \|h\|   \mbox{ for all }  h \in \tilde{C}_J(x_0),
	\end{equation*}
	This is equivalent to (\ref{robu-recov-cond}) thanks to $\tilde{C}_J(x_0)= {\rm cl}D_J(x_0)$.
	The proof of (\ref{rob-1}) is thus completed by  Lemma 2.2.

	Since $x_0$ is a sharp solution of problem $(P_J)$, there exists $y_0$ such that $A^T  y_0 \in \partial J(x_0).$ 
	Furthermore, according to \cite[Lemma 3.5]{GHS}, it can be deduced that
	\begin{equation}\label{estim-1}
		\|Ax^{\mu} -b \| \le  \delta + 2c_1\delta \|y_0\|.
	\end{equation}
	Additionally, based on \cite[Theorem 3.5.2]{IVT}, $x^{\mu}$ is proven to be a solution to the subsequent problem:
	$$\min_{x}J(x) \;\;  \mbox{ s.t. }\,\, \|Ax - b \| \le \|Ax^{\mu} - b \|.$$
	Applying the result {\rm (i)} and (\ref{estim-1}), there exists $\alpha > 0$ such that
	\begin{align*}
		\|x^{\mu}-x_0 \| \le \frac{2 \|Ax^{\mu}-b\|}{\alpha} \le \frac{2(1+2c_1 \|y_0\|)}{\alpha}\delta.
	\end{align*}
	\qed

	\bigskip
	Consider the Tikhonov regularization problem 
	$$ \min_{x}  \frac{1}{2} \|Ax - b \|^2+ \mu J(x) \leqno(QP) $$
	with a regularization parameter $\mu > 0$, and its associated Lagrange dual problem is
	$$\max_{y} - \frac{1}{2} \|y\|^2+ \langle y, b \rangle  \,\,\mbox{ s.t. } \,\,J^{o}(A^T  y) \le \mu  \Leftrightarrow A^T  y \in \mu B^o. \leqno(QD)$$
	We can easily confirm that the Lagrange dual problem of $(QD)$ is $(QP)$. By the same reasoning as in the case of $(D_J)$ and $(P_J)$ (see \cite[Theorem 3.4, 3.6]{Bonnans-Shapiro}), it can be shown that there is no duality gap between $(QD)$ and $(QP)$. Furthermore, the Lagrange multiplier set for $(QD)$ is a nonempty, convex, compact subset that coincides with the optimal solution set of $(QP)$, and is the same for any optimal solution of $(QD)$. Therefore, the uniqueness of optimal solutions for problem $(QP)$ is equivalent to the uniqueness of Lagrange multipliers for $(QD)$ for any optimal solution of $(QD)$.
	
	Optimality condition: a point $\hat{y}$ is an optimal solution of $(QD)$ if and only if there exists a Lagrange multiplier $x$, together with $\hat{y}$ satisfying
	\begin{equation*} \label{optim-condition-qd}
		Ax +\hat{y} = b \,\,\mbox{and}\,\,   x \in N_{\mu B^o}(A^T  \hat{y}),
	\end{equation*}
	or equivalently,
	\begin{equation*} \label{optim-condition-qd-1}
		Ax +\hat{y} = b \,\,\mbox{and}\,\, A^T   \frac{\hat{y}}{\mu} \in \partial J(x) = B^o \cap \{z \colon \langle x, z \rangle = J(x)\}.
	\end{equation*}
	Conversely, if $\hat{x}$ is an optimal solution to $(QP)$, then $0 \in \partial (\mu J + \frac{1}{2}\|A \cdot - b \|^2)(\hat{x}) = \mu\partial  J(\hat{x}) + A^T  (A\hat{x}-b).$ Let $\hat{y} = b-A\hat{x}$, then $A^T   \frac{\hat{y}}{\mu} \in \partial J(\hat{x}).$
	
	\begin{thm}\label{Thm 3.7} A feasible point $\hat{x}$ is the unique solution of $(QP)$ if and only if there exists $\hat{y}$ satisfying $A^T  \hat{y} \in \partial  J(\hat{x})$ such that
		\begin{equation*} \label{cond-unique-solution-qd}
			{\rm Ker}A \cap (N_{B^o}(A^T  \hat{y})+ [\hat{x}]) = \{0\}.
		\end{equation*}
	\end{thm}
	
	{\it Proof}  The proof follows a similar line of reasoning as that presented in Theorem \ref{Thm 3.1}.
	\qed
	
	\medskip 
	The equivalence of solution uniqueness for problem $(P_J)$ and $(QP)$ has been noted in previous studies, such as in \cite[Proposition 3.2]{MS} and \cite[Proposition 4.12]{FNT}. Furthermore, according to Lemma  \ref{sharp-sol}, it is evident that if $\hat{x}$ is a sharp solution of $(QP)$ with a proper lsc convex function $J$ (not necessarily a gauge), then $0 \in {\rm int} [\partial  J(\hat{x}) + A^T  (\frac{A\hat{x}-b}{\mu})] \subset {\rm int} [\partial  J(\hat{x}) + {\rm Im}A^T  ],$ which implies that $\hat{x}$ is a sharp solution of problem $(P_J)$.

	\section{Applications}
	\subsection{Sparse analysis regularization optimization}
	The convex regularization $J(\cdot )= \|D^T  \cdot \|_{1}$, where $D \colon \R^p \to \R^n$ is a dictionary, is deemed as analysis $\ell_1$-regularization in the field of inverse problems or as generalized Lasso \cite{TT} in statistics. It includes several well-known regularizers as specific instances, such as total variation \cite{ROF}, fused Lasso \cite{TSRZK}, and weighted $\ell_1$ minimization \cite{CWB}.

	In the noiseless context,   consider the following constrained problem
	$$
	\min_{x\in\R^n} \|D^T  x\|_1 \,\,\mbox{ s.t. } \,\, Ax = b_0 , \leqno(AP)
	$$
	and its Lagrange dual problem
	$$
	\max_{y\in\R^m} \langle y,b_0 \rangle \,\,\mbox{ s.t. } \,\,  J^{o}(A^T  y) \le 1,  \leqno(AD)
	$$
	where $A \colon \R^n \to \R^m$ is a linear operator and $b_0\in \R^m$ is a given vector.
	Since $\ell_1$ norm is polyhedral, by \cite[Proposition 2.2]{FMP} we have
	\begin{equation*}
		J^{o}(z) = \inf_u\{ \| u \|_{\infty} \colon Du = z \}  \mbox{ for all } z \in \R^n.
	\end{equation*}
	Furthermore, the minimum is attained when the value is finite. In this case, we have
	\begin{equation*}
		B^o = \{ z \in \R^n \colon J^{o}(z) \le 1\} = D\{ u\in \R^p \colon \|u\|_{\infty} \le 1\}.
	\end{equation*}

	Let $x_0$ be a feasible point of $(AP)$, $I^{+}=\{i \colon (D^T  x_0)_i >0 \},$  $I^{-}=\{i \colon (D^T  x_0)_i <0 \}$, and $I^0 = \{i \colon (D^T  x_0)_i =0 \}$. 
	For any $h \in \R^n,$ we have
	\begin{eqnarray*}
		\lim_{t \downarrow 0} \frac{\|D^T  (x_0 + t h) \|_1 -\|D^T  x_0\|_1}{t} & = &
		\lim_{t\downarrow 0} \frac{ \sum \limits_{i \in I^{+} \cup I^{-} } |(D^T  (x_0 + t h))_i | + \sum \limits_{i \in I^0} t| (D^T  h)_i | - \sum \limits_{i \in I^{+} \cup I^{-}} |(D^T  x_0)_i| }{t}\\
		&=&  \lim_{t\downarrow 0} \frac{ \sum \limits_{i \in I^{+} \cup I^{-}}[(D^T  (x_0 + t h))_i  - (D^T  x_0)_i] + \sum \limits_{i \in I^0} t|(D^T  h)_i|} {t}\\
		&=& \sum\limits_{i \in I^{+}}(D^T  h)_i  - \sum\limits_{i \in I^{-}}(D^T  h)_i  + \sum\limits_{i \in I^0} |(D^T  h)_i|.
	\end{eqnarray*} 	
	It can be deduced that
	$$\tilde{C}_J(x_0) = \{h \in \R^n \colon  \sum\limits_{i \in I^{+}}(D^T  h)_i + \sum\limits_{i \in I^0} |(D^T  h)_i|  \le   \sum\limits_{i \in I^{-}}(D^T  h)_i \}.$$
	As a consequence of Theorem \ref{Thm 3.5}, we have the following result. 
	\begin{cor}
		A feasible point $x_0$ is the unique solution (or a sharp solution) of problem $(AP)$
		if and only if for any $h \in {\rm Ker}A \setminus \{0\},$
		\begin{equation*}\label{s-ap-cone}
			\sum\limits_{i \in I^{+}} (D^T  h)_i + \sum\limits_{i \in I^0} |(D^T  h)_i|   > \sum\limits_{i \in I^{-}}(D^T  h)_i .
		\end{equation*}
	\end{cor}	
	
	For a given $x_0 \in \R^n$, we have 
	$$\partial J(x_0) = D\partial \|D^T  x_0\|_{1}=D\{u \in \R^p \colon u_I= {\rm sign}(D^T  x_0)_I, \|u_{I^c}\|_{\infty} \le 1\},$$
	where  $I = {\rm supp}(D^T  x_0)$ and  $I^c = \{1,2,\ldots,p\} \setminus I.$
	It is clear that
	\begin{equation*}
		{\rm ri} \partial \|D^T  x_0\|_{1} = \{u \in \R^p \colon u_I= {\rm sign}(D^T   x_0)_I, \|u_{I^c}\|_{\infty} < 1\}.
	\end{equation*}
	Let $u_0 = {\rm sign}(D^T  x_0)_I \times 0_{I^c} \in \partial \|D^T  x_0\|_1$.
	Hence, the first equality in condition {\rm (v)} in Theorem \ref{Thm 3.5}  becomes
	\begin{equation*}
		\R^n = {\rm Im}A^T   - D{\rm aff}(\partial \|D^T  x_0\|_{1}-u_0) = {\rm Im}A^T   - D(\{0_I\} \times \R^{|I^c|}),
	\end{equation*}
	which, according to polar duality, is equivalent to
	\begin{equation*}
		{\rm Ker}A \cap (D(\{0_I\} \times \R^{|I^c|}))^{\bot} =\{0\}.
	\end{equation*}
	This means that ${\rm Ker}A \cap {\rm Ker}(D^T  _{I^c})=\{0\}$. 
	Directly following Theorem  \ref{Thm 3.5}, we are able to characterize the sharp solution to $(AP)$ as follows.
	
	\begin{cor}\label{s-a-regular} 
		A feasible point  $x_0$ is the unique solution (or a sharp solution) of problem $(AP)$ if and only if
		\begin{align} \label{cond-unique-solution-8-1}
			&{\rm Ker}A \cap {\rm Ker}(D^T  _{I^c})=\{0\}, \\
			&\exists y_0 \in \R^m \,\, \mbox{and}\,\, u \in \R^p \,\,\mbox{ s.t. } \,\, A^T  y_0 =Du, u_I= {\rm sign}(D^T  {x_0})_I, \mbox{ and } \|u_{I^c}\|_{\infty} < 1.\label{cond-unique-solution-8-2}
		\end{align}
	\end{cor} 
	
	Corollary \ref{s-a-regular} presents a condition that aligns with Condition 1 from \cite{ZYY}, with the additional assumption of the surjectivity of $A$.  
	When $D$ is the identity matrix, conditions (\ref{cond-unique-solution-8-1}) and (\ref{cond-unique-solution-8-2}) simplify to those in  \cite{Gilbert,ZYC}, i.e., 
	\begin{align} \notag
		&A_I \mbox{ is injective, and}  \\
		\notag
		&\exists y_0 \in \R^m \,\,\mbox{ s.t. } \,\, A^T  _{I}y_0 = {\rm sign}(({x_0})_I)  \mbox{ and } \|A^T  _{I^c}y_0 \|_{\infty} < 1,
	\end{align}
	where $A_I$ denotes the submatrix of $A$ which is generated by the selection of columns from  $A$ based on the indices in $I={\rm supp}(x_0)$. Similar necessary and sufficient conditions have been studied in  \cite{Dossal,Fuchs,GHS,Tibshirani}.
	\medskip

	\subsection{Weighted sorted $\ell_1$-norm regularization optimization}
	The weighted sorted $\ell_1$-norm (WSL1) was introduced by Bogdan et al. \cite{BBSC} and Zeng and Figueiredo \cite{ZF} as a statistical method for promoting models with a low false discovery rate when applied to specific design matrices. It is an extension of the octagonal shrinkage and clustering algorithm for regression (OSCAR)  \cite{BR} and has good sparse clustering properties and control of the false discovery rate.
	
	The WSL1 norm with respect to the non-negative vector $w \in \R^n$, whose elements are ordered as $w_1 \ge w_2 \ge \ldots \ge w_n \ge 0$, is defined as
	\begin{equation*}
		\|x\|_w = \sum_{i=1}^n w_i|x_{[i]}|,
	\end{equation*}
	where $x_{[i]}$ is the $i$th-largest component of $x$ in magnitude. 
	Given the atomic set 
	$${\cal A} = \{ \Pi b_i \colon \Pi \in P^n, i =1,2,\ldots, n\},$$
	where $P^n$ is the set of all $n \times n$ signed permutation matrices whose elements are $0$ and $\pm 1$, and the vectors
	$$b_i := (\underbrace{\tau_i,\ldots,\tau_i}_{i \mbox{ entires}}, 0, \ldots, 0) \mbox{ with } \tau_i :=(\sum_{j=1}^iw_j)^{-1}.$$
	It follows from \cite[Theorem 1]{ZF} that
	\begin{align*}
		\|x\|_w = \gamma_{{\rm conv}\cal A}(x).
	\end{align*}
	Then ${\cal B}_w := \{x \colon \|x\|_{w} \le 1\}= {\rm conv}{\cal A}.$
	Consider the signed permutation defined as
	\begin{equation*}
		P_{w}= \{(\sigma_1 w_{\pi(1)}, \ldots, \sigma_nw_{\pi(n)})^T   \colon \sigma_1, \ldots, \sigma_n \in \{-1,1\}, \pi \in S_n\},
	\end{equation*}
	where $S_n$ denotes the set of all permutations of $\{1,2, \ldots,n\}$. We observe from \cite[Proposition 8]{STW} that ${\cal B}_w^o= ({\rm conv}{\cal A} )^o = {\rm conv}P_{w}$
	and $\partial \|x\|_{w} = \{z \in {\rm conv}P_{w} \colon \langle z,x\rangle =  \|x\|_{w}\}$.
	Let
	\begin{eqnarray*}
		{\cal E}(x) &=&\{ z \in P_{w} \colon \langle x,z \rangle = \|x\|_{w}\} \\
		&=&\Bigg \{z=(\sigma_1 w_{\pi^{-1}(1)}, \ldots, \sigma_n w_{\pi^{-1}(n)})^T   \colon
		\begin{array}{c}
			\sigma_i= {\rm sign}(x_i) \mbox{ if } x_i \ne 0;\\
			\sigma_i \in \{-1,1\} \mbox{ if } x_i =0;\\ \pi \in S_n \mbox{ satisfing } |x_{\pi(1)}| \ge \ldots \ge |x_{\pi(n)}| 
		\end{array} \Bigg\}.
	\end{eqnarray*}
	Clearly, $\gamma_{{\rm conv}P_{w}}(z)=1$ for all $z \in {\cal E}(x)$ and $\partial \|x\|_{w}={\rm conv}{\cal E}(x).$
	Hence, we obtain that
	\begin{equation*}
		{\rm aff}\partial \|x\|_{w} ={\rm aff}{\cal E}(x),
	\end{equation*}
	and
	\begin{equation*}
		{\rm ri}\partial \|x\|_{w} = \{ u \in \R^n \colon u =\sum_{z \in {\cal E}(x)} \lambda_zz, \sum_{z \in {\cal E}(x)} \lambda_z =1, \lambda_z >0  \}.
	\end{equation*}
	Consider the weighted sorted $\ell_1$-norm regularization problem 
	$$
	\min_{x \in\R^n}\|x\|_{w} \,\,\mbox{ s.t. } \,\, A x = b_0, \leqno(SP)
	$$
	and its Lagrange dual problem
	$$
	\max_{y\in\R^m} \langle y,b_0 \rangle \,\,\mbox{ s.t. } \,\,  A^T  y \in {\rm conv}P_{w}, \leqno(SD)
	$$
	where $A\colon \R^{n} \to \R^m$ is a linear mapping and $b_0 \in \R^m$. 
	According to Theorem \ref{Thm 3.5}, we can derive the following result.

	\begin{cor}
		A feasible point $\tilde{x} \in \R^n$ is a sharp solution of problem $(SP)$ if and only if 	
		\begin{align}
			&{\rm Im}A^T   - {\rm aff}({\cal E}(\tilde{x}))=\R^n, \label{cond-unique-solution-8-3-1} \\
			& \exists \tilde{y}\in \R^m \,\,\mbox{ s.t. } \,\,  A^T  \tilde{y} = \sum_{z \in {\cal E}(\tilde{x})} \lambda_zz, \sum_{z \in {\cal E}(\tilde{x})} \lambda_z =1, \mbox{and}\,\, \lambda_z >0. \nonumber
		\end{align}
	\end{cor}
	
	In particular, suppose that  $w_1 > w_2> \ldots > w_n>0$ and $\tilde{x}$ is a feasible point for problem $(SP)$. Let $J_{\tilde{x}}(i) = \{ j \in \{1,\ldots, n\} \colon |\tilde{x}_j|=|\tilde{x}_i|\}$ and $u^0 \in \R^n$ satisfy $u^0_i={\rm sign}(\tilde{x}_i)w_{\pi^{-1}(i)}$ if $\tilde{x}_i\ne 0$ and $u^0_i=w_{\pi^{-1}(i)}$ otherwise for some $\pi\in S_n$ with $|\tilde{x}_{\pi(1)}| \ge \ldots \ge |\tilde{x}_{\pi(n)}|$. Then
	\begin{eqnarray*}
		{\rm aff}({\cal E}(\tilde{x})-u^0)=  \bigg \{ v\in \R^n \colon
		\begin{array}{c}
			v_i =0 \mbox{ for all } i \mbox{ with } \tilde{x}_i \neq 0 \mbox{ and } J_{\tilde{x}}(i)=\{i\};\\
			\sum\limits_{j\in J_{\tilde{x}}(i)} {\rm sign}(\tilde{x}_j) v_j = 0 \mbox{ for all } i \mbox{ with } \tilde{x}_i \neq 0 \mbox{ and } J_{\tilde{x}}(i)\ne \{i\} 
		\end{array} \bigg\},
	\end{eqnarray*}
	and
	$$
	({\rm aff}({\cal E}(\tilde{x})-u^0))^{\perp} =\bigg\{ 
	s \in \R^n \colon 
	\begin{array}{c}
		s_i=0 \mbox{ if } \tilde{x}_i=0;\\
		s_{j}= {\rm sign}(\tilde{x}_{j})\alpha_i, j \in J_{\tilde{x}}(i), \alpha_i \in \R, \mbox{ if } \tilde{x}_i\neq 0 \mbox{ and }  J_{\tilde{x}}(i)\ne\{i\}  
	\end{array}	
	\bigg\}.
	$$
	It follows from (\ref{cond-unique-solution-8-3-1}) that 
	\begin{equation*}
		{\rm Im}A^T   - {\rm aff}({\cal E}(\tilde{x})-u^0)= 	{\rm Im}A^T   - {\rm aff}{\cal E}(\tilde{x})+u^0=\R^{n} +u^0 = \R^{n},
	\end{equation*}
	which is, by the polar duality, equivalent to ${\rm Ker}A\cap (	{\rm aff}({\cal E}(\tilde{x})-u^0))^{\perp}=\{0\}$. 
	Therefore, $\tilde{x}$ is a sharp solution of problem $(SP)$ with $w_1 > w_2> \ldots > w_n>0$ if and only if $A$ is injective restricted to  $(	{\rm aff}({\cal E}(\tilde{x})-u^0))^{\perp}$ and there exists $\tilde{y} \in \R^m$ such that
	$ A^T  \tilde{y} = \sum_{z \in {\cal E}(\tilde{x})} \lambda_zz,$ $ \sum_{z \in {\cal E}(\tilde{x})} \lambda_z =1$, and $ \lambda_z >0.$
	
	\medskip
	
	\subsection{Nuclear norm optimization}
	The nuclear norm of a matrix is the spectral analog to the vector $\ell_1$-norm, which is a gauge that promotes low rank (e.g., sparsity concerning rank-1 matrices), see \cite{RFP}. For $X\in \R^{m \times n}$, its nuclear norm is given by
	\begin{equation*}
		\|X\|_* := \sum_{i=1}^t\sigma_i(X),
	\end{equation*}
	where $\sigma_i(X) (i=1,\ldots,t)$ are the singular values of $X$ and $t = \min\{m,n\}$. 
	Hence, 
	\begin{equation*}   
		{\cal B} := \{ X \colon \|X\|_* \le 1\} \,\,\mbox{and}\,\,
		{\cal B}^o = \{Z \colon \sigma_1(Z):=\max_{1\le i \le t} \sigma_i(Z)\le 1\}.
	\end{equation*}
	Let $\Phi \colon \R^{m \times n} \to \R^p$ be a linear mapping with adjoint $ \Phi^*$ and $b_0 \in \R^p$. Consider the nuclear norm minimization problem 
	$$
	\min_{X}\|X\|_* \,\,\mbox{ s.t. } \,\,\Phi X = b_0, \leqno(NP)
	$$
	and its Lagrange dual problem
	$$
	\max_{y\in\R^p} \langle y,b_0 \rangle \,\,\mbox{ s.t. } \,\,  \Phi^* y \in {\cal B}^o. \leqno(ND)
	$$
	
	Denote by ${\cal O}^n$ the set of all $n \times n$ orthogonal matrices in $\R^{n \times n}$. Let $X_0 \in \R^{m \times n}$ with rank $r$ be a feasible solution of problem $(NP)$ and  have a singular value decomposition
	\begin{equation}\label{SVD}
		X_0 = U \Sigma(X_0) V^T  ,
	\end{equation}
	where $U=(u_1, \ldots, u_m)\in {\cal O}^m$, $V=(v_1, \ldots v_n)\in {\cal O}^n$,  and $\Sigma(X_0) \in \R^{m \times n}$ is a diagonal matrix with singular values in descending order
	$$\sigma_1(X_0) \ge \sigma_2(X_0) \ldots \ge \sigma_r(X_0) >0= \sigma_{r+1}(X_0) = \ldots = \sigma_t(X_0)$$
	on the diagonal.  Denote by ${\cal O}^{m \times n}(X_0)$ the set of such matrices $(U,V)$ in the (\ref{SVD}), i.e.,
	$${\cal O}^{m \times n }(X_0)=\{(U,V)\in {\cal O}^{m} \times {\cal O}^{n} \colon X_0 = U \Sigma(X_0) V^T  \}. $$
	
	It has been proved in \cite{lewis} that the subdifferential of $\|X_0\|_*$ is
	\begin{eqnarray*}
		\partial \|X_0\|_* &=& \{ Z = U {\rm diag}(s) V^T   \colon (U,V) \in {\cal O}^{m \times n}(X_0), s\in \partial \|\sigma(X_0)\|_1 \}\\
		&=& \{Z= U {\rm diag}(s) V^T   \colon (U,V) \in {\cal O}^{m \times n}(X_0),  s_i = 1, i=1, \ldots, r; |s_i|\le 1, i = r+1, \ldots t \},
	\end{eqnarray*}
	where ${\rm diag}(s)$ is an $m \times n$ diagonal matrix with $s$ on its diagonal. The matrices from the singular value decomposition of $X_0$ can be partitioned as
	\begin{equation}\label{partition SVD}
		U= (U^{(1)}\vdots U^{(2)})\,\,\mbox{and}\,\, V= (V^{(1)} \vdots V^{(2)}),
	\end{equation}
	where $U^{(1)}, V^{(1)}$ have $r$ columns and $U^{(2)}, V^{(2)}$ have $m-r$ and $n-r$ columns, respectively. It was shown  by \cite[Example 2]{watson} that the subdifferential of $\|X_0\|_{*}$ can be represented as
	\begin{equation}\label{subdif-nun}
		\partial \|X_0\|_* = \{Z = U^{(1)}(V^{(1)})^T  + U^{(2)}W^{(2)}(V^{(2)})^T   \colon ((U^{(1)}\vdots U^{(2)}),(V^{(1)} \vdots V^{(2)})) \in {\cal O}^{m \times n}(X_0), \sigma_1(W^{(2)}) \le 1\},
	\end{equation}
	where $W^{(2)}$ is $(m-r) \times (n-r)$ matrix.
	Observe that $U^{(1)}(V^{(1)})^T $ is invariant for any singular value decomposition of $X_0$. Indeed, suppose that this matrix $X_0$ has another singular value decomposition $X_0 = {\overline U} \Sigma(X_0) {\overline V^T  }$ with $(\overline U,\overline V) \in {\cal O}^{m \times n}(X_0)$. Since the singular values are distinct, the vector $\sigma(X_0)$ can have the following form:
	\begin{eqnarray*}
		& &[\sigma_1(X_0), \sigma_2(X_0), \ldots, \sigma_t(X_0)]^T  \\
		&=&[\underbrace{\mu_1(X_0),\ldots, \mu_1(X_0)}_{r_1}, \underbrace{\mu_2(X_0),\ldots, \mu_2(X_0)}_{r_2},\ldots, \underbrace{\mu_q(X_0),\ldots, \mu_q(X_0)}_{r_q},\underbrace{0,\ldots, 0}_{t-r}]^T  ,
	\end{eqnarray*}
	where
	$$
	\mu_1(X_0)> \mu_2(X_0)>\ldots> \mu_q(X_0)>0, \,\, r_1+r_2+\ldots+r_q=r,
	$$
	and $r_i$ is the number of $\mu_i(X_0)$ for $i=1,\ldots,q$.
	It follows from singular value decomposition that $X_0V_i=\mu_i(X_0)U_i$ and $X_0\overline{V}_i=\mu_i(X_0)\overline{U}_i$,
	where $U_i, \overline{U}_i\in \mathbb{R}^{m\times r_i}$ and $V_i, \overline{V}_i\in \mathbb{R}^{n\times r_i}$ have orthonormal columns for any $i=1,\ldots,q$.
	Then there exist orthonormal matrices $T_i$ and $S_i$ of order $r_i$ such that $\overline{U}_i=U_iT_i$ and  $\overline{V}_i=V_iS_i$.
	It holds that
	$$\mu_i(X_0) U_iT_i=\mu_i(X_0)\overline{U}_i=X_0\overline{V}_i=X_0V_iS_i=\mu_i(X_0) U_iS_i \,\,\mbox{for any}\,\, i=1,\ldots,q.$$
	Since the column vectors of $U_i$ are mutually orthogonal and $\mu_i(X_0)>0$, we have $T_i=S_i$ for any $i=1,\ldots,q$.
	It follows that
	\begin{equation*}
		\overline{U}^{(1)}	(\overline{V}^{(1)})^T   = \sum_{i=1}^q  U_iT_{i} (V_iS_{i})^T  
		= \sum_{i=1}^q U_iT_{i}S_{i}^T   V_i^T  
		= \sum_{i=1}^q U_i V_i^T  = U^{(1)}	(V^{(1)})^T .
	\end{equation*}
	Thus given any singular value decomposition of $X_0 = (U^{(1)}\vdots U^{(2)})\Sigma(X_0) (V^{(1)} \vdots V^{(2)})^T  $, (\ref{subdif-nun}) can be equivalently written as 
	\begin{equation} \label{watson-sub}
		\partial \|X_0\|_* = \{Z = U^{(1)}(V^{(1)})^T  + U^{(2)}W^{(2)}(V^{(2)})^T  \colon  W^{(2)} \in \R^{(m-r) \times (n-r)}, \sigma_1(W^{(2)}) \le 1\}.
	\end{equation}
	
	It is known that ${\cal B}^o = \partial \| \mathsf{0}\|_* = \{ Z \in \R^{m \times n} \colon \sigma_1(Z) \le 1 \}.$
	By \cite[Corollary 2.5]{lewis}, the normal cone to ${\cal B}^o$ at a matrix $Z$ with $\sigma_1(Z)=1$ is
	\begin{equation*}
		N_{{\cal B}^o}(Z) = \partial (I_{B_{\infty}}\circ \sigma)(Z) = \{U {\rm diag}(s) V^T   \colon  (U,V)\in {\cal O}^{m\times n}(Z), s \in N_{B_{\infty}}(\sigma(Z))\}.
	\end{equation*}
	Denote by $p(Z) :=\#\{i  \colon \sigma_i(Z) = 1\}$ the number of $i$ with $\sigma_i(Z)=1$. Since $N_{[-1,1]}(1) = \R_{+}$ and $N_{[-1,1]}(\rho)=\{0\}$ for any $0\le \rho  <1$. Hence, we obtain that  $s \in N_{B_{\infty}}(\sigma(Z))$ if and only if
	$s_i \in N_{[-1,1]}(\sigma_i(Z))$, $i=1, \ldots,t$.
	
	Suppose that there exists $y_0\in \R^p$ such that $\Phi^* y_0 \in \partial \|X_0\|_{*}$.  It follows from \cite[Corollary 2.5]{lewis} that $\sigma(\Phi^* y_0) \in \partial \| \sigma(X_0)\|_1$ and there exists a simultaneous singular value decomposition of the form
	\begin{equation*}
		X_0 = U \Sigma(X_0) V^T   \,\,\mbox{and}\,\,  \Phi^* y_0 = U \Sigma(\Phi^* y_0) V^T    
	\end{equation*}
	with $U \in {\cal O}^m$ and $V \in {\cal O}^n$. 
	Then by \cite[Proposition 10]{ZS} we have
	$N_{{\cal B}^o}(\Phi^* y_0) = U 
	\begin{pmatrix}
		S_{+}^{p(\Phi^*y_0)} & 0 \\
		0 & 0\\
	\end{pmatrix} V^T  $,
	where $S_{+}^p :=\{ X \in {\cal S}^p \colon X \succeq 0 \}$ denotes the set of all $p \times p$ symmetric positive semidefinite matrices and ${\cal S}^p$ denotes the set of all $p \times p$ symmetric matrices.   
	According to Theorem \ref{Thm 3.1}, we can characterize the unique solution to $(NP)$ as follows.
	\begin{cor}  A feasible point $X_0$ is the unique minimizer of $(NP)$ if and only if there exists $y_0$ satisfying $\Phi^*y_0 \in \partial \|X_0\|_*$ such that	
		\begin{equation*} \label{radial-unique}
			{\rm Ker} \Phi \cap  \Bigg( U 
			\begin{pmatrix}
				S_{+}^{p(\Phi^*y_0)} & 0 \\
				0 & 0\\
			\end{pmatrix} V^T   +[X_0] \Bigg)  = \{0\}.
		\end{equation*}
	\end{cor}  
	\medskip

	The assertion has been proved by \cite[Lemma 3.1]{LPB} and \cite[Corollary 4.8]{FNP-1}, and an explicit formula of the radial cone ${\cal R}_{N_{{\cal B}^o}(\Phi^* y_0)}(X_0) $ has been presented in \cite[Lemma 3.5]{FNP-1}. Additionally, \cite[Theorem 3.1, 3.2]{LPB} also offers some sufficient conditions for the unique solution.

	\cite[Theorem 5.2]{FNP} illustrates that 
	$X_0$ is a strong solution to $(NP)$, i.e., there exists $c >0$ and $\delta >0$ such that 
	$$\|X\|_* - \|X_0\|_* \ge c \|X -X_0\|^2_F \mbox{ for all }  \|X - X_0\|_F < \delta \mbox{ and } \Phi X =b_0,$$ 
	if and only if 
	$ {\rm Ker} \Phi \cap {\rm cl} (N_{{\cal B}^o}(\Phi^* y_0) + [X_0]) = {\rm Ker} \Phi \cap T_{N_{{\cal B}^o}(\Phi^* y_0)}(X_0) =\{0\}$ and
	\cite[Corollary 4.2]{FNP} gives an explicit formula of $T_{N_{{\cal B}^o}(\Phi^* y_0)}(X_0)$.
	Since the subdifferential mapping of the nuclear norm is metrically subregular, we can conclude that $X_0$ is a strong solution to $(NP)$ if and only if any of the equivalent conditions from Theorem \ref{Thm 3.2} and Theorem \ref{Thm 3.3} holds.
	
	As a direct result of Theorem \ref{Thm 3.5}, we get the following characterization of a sharp solution of $(NP)$.
	\begin{cor} 
		Let $X_0 \in \R^{m \times n}$   with rank $r$ be a feasible solution of problem $(NP)$ and have a  singular value decomposition (\ref{SVD}) with $U,V$ in (\ref{partition SVD}). Then
		$X_0$  is a sharp solution of $(NP)$ if and only if there exists $y_0$ satisfying $\Phi^*y_0 \in \partial \|X_0\|_*$ such that any of the
		following equivalent conditions holds:\\
		{\rm(i)} 
		\begin{equation*}
			{\rm Ker}\Phi \cap  U^{(2)} \tilde{U} 
			\begin{pmatrix}
				S_{+}^{p((U^{(2)})^T  (\Phi^* y_0) V^{(2)})} & 0 \\
				0 & 0\\
			\end{pmatrix} \tilde{V}^T   (V^{(2)})^T   =\{0\},
		\end{equation*}
		where $(\tilde{U},\tilde{V}) \in {\cal O}^{(m-r) \times (n-r)}((U^{(2)})^T  (\Phi^* y_0) V^{(2)});$ \\
		{\rm (ii)} $\Phi$ is injective restricted to the subspace  $T:=\{Y \in \R^{m \times n} \colon (U^{(2)})^T   Y V^{(2)} =0\}$ and there exists $z\in \R^p$ and $W^{(2)} \in \R^{(m-r) \times (n-r)}$  with  $\sigma_1(W^{(2)}) < 1$  such that $\Phi^* z =U^{(1)}(V^{(1)})^T  + U^{(2)}W^{(2)}(V^{(2)})^T $;\\
		{\rm (iii)} For any $H \in {\rm Ker}\Phi \setminus \{0\}$, we have 
		\begin{equation*}
			{\rm tr}((U^{(1)})^T   H V^{(1)}) + \|(U^{(2)})^T  H V^{(2)}\|_* > 0.
		\end{equation*} 	 
	\end{cor}
	{\it Proof}	 It follows from (\ref{watson-sub}) that there exists $W_0^{(2)} \in \R^{(m-r) \times (n-r)}$ with $\sigma_1(W_0^{(2)}) \le 1$ such that
	$\Phi^* y_0 =U^{(1)}(V^{(1)})^T  + U^{(2)}W_0^{(2)}(V^{(2)})^T $. This implies that $W_0^{(2)} = (U^{(2)})^T  (\Phi^* y_0) V^{(2)}.$
	It is easy to verify that
	\begin{equation*}
		N_{\partial \| X_0\|_{*}}(\Phi^* y_0) =\{ H \in \R^{m \times n} \colon (U^{(2)})^T  H V^{(2)} \in N_{ {{\cal B}^{(2)}}^o} ((U^{(2)})^T  (\Phi^* y_0) V^{(2)})\},
	\end{equation*}
	where
	$ {{\cal B}^{(2)}}^o = \{W^{(2)} \in \R^{(m-r) \times (n-r)}  \colon \sigma_1(W^{(2)}) \le 1\}.$ 
	Let $$(U^{(2)})^T  (\Phi^* y_0) V^{(2)}= \tilde{U} \Sigma((U^{(2)})^T  (\Phi^* y_0) V^{(2)}) \tilde{V}^T  $$ be a singular value decomposition, where
	$\tilde{U}  \in {\cal O}^{m-r} $ and $ \tilde{V} \in {\cal O}^ {n-r}$. 
	Then we obtain that 
	\begin{equation*}
		N_{ {{\cal B}^{(2)}}^o} ((U^{(2)})^T  (\Phi^* y_0) V^{(2)}) =  \tilde{U} 
		\begin{pmatrix}
			S_{+}^{p((U^{(2)})^T  (\Phi^* y_0) V^{(2)}) } & 0 \\
			0 & 0\\
		\end{pmatrix} \tilde{V}^T  
	\end{equation*}
	and 
	\begin{equation*}
		N_{\partial \| X_0\|_{*}}(\Phi^* y_0)= U^{(2)} \tilde{U} 
		\begin{pmatrix}
			S_{+}^{p((U^{(2)})^T  (\Phi^* y_0) V^{(2)}) } & 0 \\
			0 & 0\\
		\end{pmatrix} \tilde{V}^T   ( V^{(2)})^T . 
	\end{equation*}
	Therefore, {\rm (i)} follows from Theorem \ref{Thm 3.5}{\rm (i)}. 	
	
	It is easy to see that
	\begin{eqnarray*}
		{\rm aff} \partial \|X_0\|_* =  \{Z = U^{(1)}(V^{(1)})^T  + U^{(2)}W^{(2)}(V^{(2)})^T  \colon  W^{(2)} \in \R^{(m-r) \times (n-r)}\}.
	\end{eqnarray*}
	It is clear that  $Z^0 := U^{(1)}(V^{(1)})^T   \in \partial \|X_0\|_*$. Hence, we have
	\begin{equation*}
		\R^{m \times n} = {\rm Im}\Phi^* - {\rm aff}(\partial \|X_0\|_{*}-Z^0) = {\rm Im} \Phi^* - \{ U^{(2)}W^{(2)}(V^{(2)})^T  \colon W^{(2)} \in \R^{(m-r) \times (n-r)} \},
	\end{equation*}
	which, according to polar duality, is equivalent to
	\begin{equation*}
		{\rm Ker}\Phi \cap \{U^{(2)}W^{(2)}(V^{(2)})^T  \colon W^{(2)} \in \R^{(m-r) \times (n-r)} \}^{\bot} =\{0\}.
	\end{equation*}
	That is ${\rm Ker}\Phi \cap \{ Y \in \R^{m \times n} \colon (U^{(2)})^T  YV^{(2)} =0\}=\{0 \}$, i.e., $\Phi$ is injective restricted to the subspace  $T$.
	Let $F =\{ s \in \R^t \colon s_i =1, i =1, 2, \ldots,r; |s_i| \le 1, i = r+1, \ldots, t\}$. Then by \cite[Example 5.7]{sa}, $F$ is a standard face of closed unit ball $B_{\infty}$ of $\R^t$ with the vector norm $\|\cdot\|_{\infty}$ associated with $\partial \|X_0\|_*$. Furthermore, \cite[Theorem 6.2]{sa} says that
	\begin{eqnarray*}
		{\rm ri} \partial \|X_0\|_* &=& \{ Z \in \partial \|X_0\|_* \colon \sigma (Z) \in {\rm ri}F\}\\
		&=& \{Z = U^{(1)}(V^{(1)})^T   + U^{(2)}W^{(2)}(V^{(2)})^T  \colon W^{(2)} \in \R^{(m-r) \times (n-r)}, \sigma_1(W^{(2)}) < 1\}.
	\end{eqnarray*}
	Thus {\rm (ii)} follows from Theorem \ref{Thm 3.5}{\rm (v)}.  
	
	According to \cite[Theorem 3.1]{ZZX}, it is easy to see that
	\begin{equation*}
		\|\cdot\|_*^{\prime}(X_0,H) = {\rm tr}((U^{(1)})^T   H V^{(1)}) + \|(U^{(2)})^T  HV^{(2)}\|_*.
	\end{equation*}  
	Hence,  {\rm (iii)} follows from Theorem \ref{Thm 3.5}{\rm (vii)}. 
	\qed  
	
	\medskip 
	
	Condition {\rm (ii)} is exactly conditions 1 and 2 in \cite[Lemma 3.1]{CR}. Additionally, \cite[Theorem 4.6]{FNT} provides equivalent characterizations of sharp minima.
	
	As a consequence of Theorem \ref{Thm 3.5}, we obtain the following sufficient condition for a sharp solution of $(NP)$.

	\begin{cor} 	Let $X_0 \in \R^{m \times n}$   with rank $r$ be a feasible solution of problem $(NP)$. 
		If 
		\begin{equation} \label{a}
			\sum_{i=1}^{r}\sigma_i(H) < \sum_{i=r+1}^{t}\sigma_i(H)
		\end{equation}
		for all $H \in {\rm Ker}\Phi \setminus \{0\}$, 
		then $X_0$ is a sharp minimum of $(NP)$.
	\end{cor}
	{\it Proof}  We claim that 
	$$C_{\|\cdot\|_*}(X_0) \subset \{ H \in \R^{m \times n} \colon \sum_{i=1}^{r}\sigma_i(H) \ge \sum_{i=r+1}^t\sigma_i(H) \}.$$
	Indeed, it follows from \cite[Lemma A.18]{FR} that for any $\tau \downarrow 0$, we have
	\begin{eqnarray*}
		\|X_0 + \tau H\|_* - \|X_0\|_* &=&\sum_{i=1}^{t}\sigma_i(X_0 - \tau (-H)) -\sum_{i=1}^{r}\sigma_i(X_0)\\ 
		&\ge & \sum_{i=1}^{t}|\sigma_i(X_0) - \tau\sigma_i(H)|  -\sum_{i=1}^{r}\sigma_i(X_0)\\
		&= & \sum_{i=1}^{r}|\sigma_i(X_0) - \tau\sigma_i(H)| +  \sum_{i=r+1}^{t}|\sigma_i(X_0) - \tau\sigma_i(H)|  -\sum_{i=1}^{r}\sigma_i(X_0)\\
		&\ge & -\tau\sum_{i=1}^{r}\sigma_i(H) +  \tau \sum_{i=r+1}^{t}\sigma_i(H),
	\end{eqnarray*}
	and then
	\begin{equation*}
		-\sum_{i=1}^{r}\sigma_i(H) + \sum_{i=r+1}^{t}\sigma_i(H) \le \lim_{\tau \downarrow 0} \frac{\|X_0 + \tau H\|_* - \|X_0\|_*}{\tau}. 
	\end{equation*}
	According to condition {\rm (vii)} from Theorem \ref{Thm 3.5}, the proof is complete.
	\qed 
	
	\cite[Theorem 4.40]{FR} shows that condition (\ref{a}) is equivalent to every feasible matrix $X$ of rank at most $r$ being the unique solution of problem $(NP)$.

	\medskip
	\subsection{Conic gauge optimization}
	We consider a gauge regularization optimization problem with a conic constraint as follows:
	$$
	\min_{x \in K} \tilde{J}(x) \;\; \mbox{ s.t. } \;\; Ax = b_0, \leqno(CP)
	$$
	where $\tilde{J}$ is a gauge, $K\subset \R^n$ is a closed convex cone, $A\colon \R^{n} \to \R^m$ is a linear mapping, and $b_0 \in \R^m$.
	This problem can be reformulated as a gauge optimization problem $(P_J)$ by defining $J = \tilde{J} + I_K$.
	It is clear that $J$ is a gauge and is polyhedral if $\tilde{J}$ and $K$ are polyhedral.
	By \cite[Corollary 2.5]{FMP},  if either ${\rm ri (dom }\tilde{J}) \cap {\rm ri}K \ne \emptyset$ or $\tilde{J}$ and $K$ are polyhedral, then
	\begin{equation*}
		J^{o}(z) = \inf_{u}\{\tilde{J}^{o}(z-u)+ I_{K^{o}}(u)\},
	\end{equation*}
	and moreover, the infimum is attained when finite. If $J(x_0)$ is finite, ${\rm ri (dom }\tilde{J}) \cap {\rm ri}K \ne \emptyset$ or $\tilde{J}$,  $K$ are polyhedral and ${\rm dom}\tilde{J} \cap K \ne \emptyset$, we have that $\partial J(x_0) = \partial \tilde{J}(x_0) + N_K(x_0)$ and ${\rm d}J(x_0)(h) = {\rm d} \tilde{J}(x_0)(h) + I_{T_K(x_0)}(h)$ by \cite[Corollary 10.9]{RW}.
	Hence, the critical cone of $J$ at $x_0$ is $\tilde{C}_J(x_0) = \{h \in T_K(x_0)\colon {\rm d}\tilde{J}(x_0)(h) \le 0 \}.$
	The Lagrange dual problem becomes
	$$
	\max_{y\in\R^m} \langle y,b_0 \rangle  \;\; \mbox{ s.t. } \;\;  A^T  y \in B^o +K^o, \leqno(CD)
	$$
	where $B^o = \{ z \colon \tilde{J}^o(z) \le 1\}. $
	
	The following result gives the characterization of sharp solution of problem $(CP)$ with $\tilde{J}=\|\cdot\|_1$ and $K=\R^n_{+}$.
	\begin{cor}
		A solution $x_0$ is a sharp solution of the following problem 
		\begin{equation}\label{l1-cone-p} 
			\min_{x \in \R^n_{+}} \|x\|_1  \;\; \mbox{ s.t. } \;\;  Ax = b_0,
		\end{equation}
		if and only if,  for any $h \in {\rm Ker}A \setminus \{0\},$ 
		either $h_i <0 $ for  some $i \in I^c$ or $h_{I^c} \ge 0 \Rightarrow \sum\limits_{i \in I}h_i + \sum\limits_{i \in I^c} h_i >0,$
		where $I = {\rm supp}(x_0)$ and $I^c$ is its complementary.
	\end{cor}	
	{\it Proof}	 For any $h \in \R^n,$ we have
	\begin{eqnarray*}
		\lim_{t \downarrow 0} \frac{\|x_0 + t h \|_1 -\|x_0\|_1}{t}
		= \sum\limits_{i \in I}h_i + \sum\limits_{i \in I^c} |h_i|.
	\end{eqnarray*} 	
	It can be deduced that
	$$\tilde{C}_J(x_0) = \{h \in \R^n \colon h_{I^c} \ge 0,  \sum\limits_{i \in I}h_i + \sum\limits_{i \in I^c} |h_i|  \le  0 \}.$$
	Therefore, the conclusion follows from Theorem \ref{Thm 3.5}.
	\qed 
	
	\cite[Exercise 4.8]{FR} shows that  every nonnegative $s$-sparse vector $x$ is the unique solution of problem (\ref{l1-cone-p}) 
	if and only if $h_{I^c} \ge 0 \Rightarrow \sum_{i=1}^nh_i >0$ for all $h \in {\rm Ker}A \setminus \{0\}$ and all $I \subset \{1, 2, \ldots,n\}$ with ${\rm card}(I) \le s .$

	\subsection{Semidefinite conic gauge optimization}
	Consider the semidefinite optimization problem
	$$
	\min_{X \in {\cal S}^n_{+}} \langle C, X \rangle \;\; \mbox{ s.t. } \;\;  \Phi X = b_0, \leqno(SP)
	$$ 
	where $C\in {\cal S}^n_{+} $, $\Phi \colon {\cal S}^n \to \R^m$ is a linear map, and $b_0 \in \R^m$.
	Define $J(X) = \langle C, X \rangle\ + I_{{\cal S}^n_{+}}(X)$ for any $X\in \R^{n \times n}$. Then $J$ is a gauge and its polar is
	\begin{equation*}
		J^{o}(Z) = \inf\{ \alpha \ge 0 \colon \alpha C -Z \in {\cal S}^n_{+} \}
		= \max\{0, \lambda_{\rm max}(Z,C)\},
	\end{equation*}
	where $\lambda_{\rm max}(Z,C)$ is the largest generalized eigenvalue corresponding to the eigenvalue problem $Zx= \lambda Cx,$ see \cite{FMP}. In particular, whenever $C= I$, the problem $(SP)$ is the problem of minimizing the trace of a positive semidefinite matrix with a linear constraint. It is clear that the critical cone of $J$ at $X_0$ is 
	\begin{equation*}
		\tilde{C}_J(X_0) = \{H \in T_{{\cal S}^n_{+}}(X_0)\colon  \langle C, H \rangle \le 0 \} =\{H \in {\cal S}^n \colon   E^T  HE \succeq 0 ,\langle C, H \rangle \le 0\}, 
	\end{equation*}
	where $E$ is an $n \times (n-r)$ matrix of full columns rank $n-r$ with $r = {\rm rank}(X)$ such that $XE =0$, see \cite{Shapiro}.
	As a consequence of Theorem \ref{Thm 3.5}, we get the following result.
	\begin{cor} 	A feasible point $X_0$ is a sharp minimum of problem $(SP)$ 			
		if and only if for any $H  \in {\cal S}^n \cap {\rm Ker}\Phi \setminus \{0\}$,  either $ E^T  HE \not \succeq 0$ or $\langle C, H \rangle  > 0.$	
	\end{cor}

	\section{Conclusion}	
	This paper examines the uniqueness, sharpness, and robust recovery of solutions in sparse regularization $(P_J)$ with a closed gauge from a dual point of view. We note that the Lagrange dual problem of the dual problem $(D_J)$ is the original problem $(P_J)$ and the strong duality holds. By exploiting the uniqueness of the Lagrange multipliers for the dual problem, we give a characterization of the unique solution to $(P_J)$ via the so-called radial cone. Furthermore, we present characterizations of the unique solution together with the calmness properties of two types of solution mappings. We also prove the equivalence of these characterizations under the assumption of metric subregularity of the subdifferentials. In addition, we give some sufficient and necessary conditions for a sharp solution, which ensures robust recovery with a linear rate and implies local upper Lipschitzian properties of the solution mapping. Specifically, the strict Robinson's condition for problem $(D_J)$, along with its simple reformulation (condition ({\rm vii})), has been extensively studied in the literature, particularly in the context of polyhedral gauges such as $\ell_1$ norm, the $\ell_1$-analysis model, the nuclear norm, and non-negative continuous convex function. We not only extend this condition to the general case whenever the regularizer is a non-negative lsc convex function, but also get several equivalent conditions. Applications to various regularization methods are explored in the last section.
	
	\textrm{\textbf{Acknowledgements}}
		Authors would like to thank professor Pan Shaohua and Liu Yulan for their carefully reading and valuable discussion, which improved the first version of this paper. We also thank professor Tran T. A. Nghia. for sending the preprint of Ref. 18.


\end{document}